\documentclass[psamsfonts,11pt]{amsart}

\newtheorem{theorem}{Theorem}[section]
\newtheorem{lemma}[theorem]{Lemma}
\newtheorem{proposition}[theorem]{Proposition}
\newtheorem{corollary}[theorem]{Corollary}

\theoremstyle{definition}
\newtheorem{definition}[theorem]{Definition}
\newtheorem{example}[theorem]{Example}

\newtheorem{remark}[theorem]{Remark}

\newtheorem{situation}[theorem]{Situation}

\theoremstyle{remark}

\numberwithin{equation}{section}

\newcommand{\ZZ}{\mathbb{Z}}

\newcommand{\PP}{\mathbb{P}}
\renewcommand{\AA}{\mathbb{A}}

\newcommand  {\shF}     {\mathcal{F}}
\newcommand  {\shG}     {\mathcal{G}}

\newcommand  {\shM}     {\mathcal{M}}

\newcommand  {\shN}     {\mathcal{N}}
\newcommand  {\shL}     {\mathcal{L}}
\newcommand  {\shR}     {\mathcal{R}}
\newcommand  {\shS}     {\mathcal{S}}
\newcommand  {\shT}     {\mathcal{T}}

\newcommand  {\shQ}     {\mathcal{Q}}

\newcommand  {\Char}    {\operatorname{char}}

\newcommand  {\dual}    {\vee}

\newcommand  {\im}      {\operatorname{im}}

\newcommand  {\lra}     {\longrightarrow}

\renewcommand{\O}       {\mathcal{O}}

\newcommand  {\Proj}    {\operatorname{Proj}}

\newcommand  {\ra}      {\rightarrow}

\newcommand  {\rank}    {\operatorname{rank}}

\newcommand  {\Spec}    {\operatorname{Spec}}

\newcommand  {\Syz}     {\operatorname{Syz}}

\newcommand {\soclo}{\star}
\newcommand {\pasoclo}{\star}

\newcommand {\mindeg}{\rho}
\newcommand {\maxdeg}{\lambda}
\newcommand {\gr} {+ \rm gr}

\def\mydate{\number\day\space\ifcase\month \or January\or February\or March\or April\or May\or
June\or July\or August\or September\or October\or November\or
December\fi \space\number\year}

\usepackage{amscd}
\usepackage{amssymb}

\begin{document}

\setlength{\parindent}{0cm}

\title{Computing the tight closure in dimension two}


\author[Holger Brenner]{Holger Brenner}
\address{Department of Pure Mathematics, University of Sheffield,
  Hicks Building, Houns\-field Road, Sheffield S3 7RH, United Kingdom}
\email{H.Brenner@sheffield.ac.uk}





\begin{abstract}
We study computational aspects of the tight closure of a homogeneous
primary ideal in  a two-dimensional normal standard-graded domain.
We show how to use slope criteria for the sheaf of syzygies for
generators of the ideal to compute its tight closure. In particular,
our method gives  an algorithm to compute the tight closure of three
elements under the condition that we are able to compute the
Harder-Narasimhan filtration. We apply this to the computation of
$(x^{a},y^{a},z^{a})^*$ in $K[x,y,z]/(F)$, where $F$ is a
homogeneous polynomial.
\end{abstract}

\maketitle


\section*{Introduction}

\medskip
Let $I \subseteq R$ denote an ideal in a Noetherian domain $R$
over a field $K$ of characteristic $p >0$.
The tight closure of $I=(f_1, \ldots ,f_n)$ is again an ideal
defined by
$$I^* = \{ f\in R : \exists c \neq 0: c f^{q} \in (f_1^q, \ldots ,f_n^{q})
\mbox{ holds for almost all } q=p^{e} \} \, .$$ The theory of tight
closure was developed by Hochster and Huneke (see
\cite{hochsterhunekebriancon}, \cite{hunekeapplication},
\cite{hunekeparameter}, \cite{smithtightintroduction}) and has many
applications in commutative algebra, homological algebra and
algebraic geometry. Its strength lies in the interplay of inclusion
and exclusion results for tight closure.

Huneke writes, ``Tight closure is very difficult to compute; indeed
that is necessarily the case. It contains a great deal of
information concerning subtle properties of the ring and the ideal''
(\cite[Basic Notions]{hunekeapplication}). The problem lies in the
fact that due to the definition we have to check infinitely many
conditions.

If the ring $R$ is regular, then $I=I^*$ holds for every ideal $I
\subseteq R$. If the ring is one-dimensional, then $I^* =R \cap
IR^{\rm nor}$ (the normalization), so in these two cases the
computation of $I^*$ is easy (at least there is a translation to
other more elementary computational problems). But even in the case
of a normal two-dimensional standard-graded domain $R$ over an
algebraically closed field $K$ (e. g. $R=K[x,y,z]/(F)$, $F$
homogeneous) very little is known. The tight closure of a
homogeneous parameter ideal $(f_1,f_2) \subset R$ is given by
$$ (f_1,f_2)^*= (f_1,f_2) + R_{\geq \deg(f_1) + \deg (f_2)} \, $$
(in characteristic $0$ or $p \gg 0$). This is the Strong Vanishing
Theorem of Huneke and Smith (\cite{hunekesmithkodaira}), which has
been generalized for parameter ideals in higher dimensions by Hara
(\cite{hararationalfrobenius}).

For a homogeneous $R_+$-primary ideal $I=(f_1, \ldots ,f_n)
\subseteq R$ not much is known about the tight closure $(f_1,
\ldots, f_n)^*$ for $n \geq 3$. A difficult but elementary
computation due to Singh shows that $xyz \in (x^2,y^2,z^2)^*$ holds
in the Fermat cubic given by $x^3+y^3+z^3=0$
(\cite{singhcomputation}, \cite{leuschkeappendix}). Smith has given
the two degree bounds $R_{\geq  2 (\max_i \{ \deg (f_i) \}) }
\subseteq I^*$ and $R_{\geq \deg (f_1)+ \ldots + \deg (f_n)}
\subseteq I^*$ (\cite{smithgraded}). Smith also proved a degree
bound from below: if $\deg (f_0) \leq \min \{\deg (f_1), \ldots,
\deg (f_n)\}$, then $f_0 \in I^*$ if and only if $f_0 \in I$
already. All these bounds are rather coarse for non-parameter
ideals. If the $f_i$ have the same degree $d$, then these degree
bounds say nothing between $d$ and $2d$. In particular they do not
yield anything interesting in the example of Singh.

\medskip
Another approach was initiated by Katzman and further developed by
Sullivant. They use an algorithm which computes, given an ideal $I$,
two ideals $I_1$ and $I_2$ such that $I_1 \subseteq I^* \subseteq
I_2$. If both approximations coincide, then the algorithm gives the
right answer for a fixed prime number $p$. The computations of
Sullivant (implemented in Macaulay2) of $(x^3,y^3,z^3)^*$ for the
Fermat rings $K[x,y,z]/(x^d+y^d+z^d)$ for $p \leq 53$ and $d \leq
26$ are striking and have led to some interesting observations and
conjectures. Of course one cannot expect any general result by this
method, and the lower bound $I_1$ computes rather the Frobenius
closure of the ideal. Therefore the algorithm does not give the
right answer for $(x^2,y^2,z^2)^*$ for $p = 1 \, {\rm mod }\, 3$.

\medskip
In this paper we want to
attack the problem of computing the tight closure of an ideal
from another point of view:
using the slope criteria for vector bundles.
This rests upon the geometric interpretation of tight closure via
vector bundles and projective bundles which we have developed
in \cite{brennertightproj} and \cite{brennerslope}.
This paper will emphasize the computational usefulness of this approach.

The main object to consider in this approach is the sheaf of
syzygies of total degree $m$ for homogeneous ideal generators $f_1,
\ldots ,f_n$. This is a locally free sheaf $\Syz(m)$ on the smooth
projective curve $Y = \Proj\, R$. The slope properties of this sheaf
are crucial for the underlying tight closure problem. So we may
forget the definition of tight closure and struggle instead with the
notions of slope, minimal and maximal slope, semistability and the
Harder-Narasimhan filtration of this sheaf of syzygies. This is
still a difficult task, however we can use many more tools from
algebraic geometry to attack the tight closure problem. We recall
this geometric interpretation and the resulting slope criteria for
tight closure briefly in Section \ref{sectionrecall}.

In fact we work with the notion of solid closure (denoted
$I^\star$), which coincides with tight closure in positive
characteristic and gives a satisfactory notion for characteristic
zero in dimension two. The slope conditions are easier to formulate
in zero characteristic; hence we restrict largely to this case in
the introduction.

If the sheaf of syzygies $\Syz(m)$ is semistable, then we have an
easy numerical criterion for tight closure: the common degree bound
for inclusion and exclusion is given by $\frac{\deg (f_1) + \ldots +
\deg (f_n)}{n-1}$. If the sheaf of syzygies is not semistable, then
we can argue along the Harder-Narasimhan filtration of $\Syz(m)$.
This yields the first step of an algorithm to compute the tight
closure. This algorithm gives a complete answer if the maximal
destabilizing subsheaf of the sheaf of syzygies is itself
semistable. This condition is of course fulfilled if the rank of
$\Syz(m)$ is two; hence we get an algorithm to compute the tight
closure of a homogeneous $R_+$-primary ideal generated by three
elements (section \ref{algorithmsection}) --- at least if we are
able to compute the Harder-Narasimhan filtration, which means for
rank two to find the invertible subsheaves of maximal degree
(equivalently, to compute the $e$-invariant of the corresponding
ruled surface). This method works also in positive characteristic
and we get (less complete) results about the plus closure (Section
\ref{tightthree}).

In Section \ref{tightrelations} we study the global sections of the
sheaf of syzygies for homogeneous ideal generators. Their existence
and nonexistence in certain degrees has many consequences on the
slope properties of $\Syz(m)$ and therefore on the tight closure.
For example, if there does not exist a syzygy $\neq 0$ of total
degree $k$ for the homogeneous primary elements $f_1,f_2,f_3$, then
$R_m \subseteq (f_1,f_2,f_3)^\star$ holds for $m \geq d_1+d_2+d_3 -k
+ \frac{g-1}{\delta}$, where $g$ denotes the genus and $\delta $ the
degree of $Y$ (Corollary \ref{relationboundinclusion}). The
nonexistence of global syzygies also implies exclusion results: if
there does not exist a non-trivial global syzygy for $f_1, \ldots
,f_n \in R=K[x,y,z]/(F)$ of degree $k$, then $(f_1, \ldots ,
f_n)^\star \cap R_m = (f_1, \ldots, f_n) \cap R_m$ holds for $m \leq
k- \delta +2$, where $\delta $ is the degree of $F$ (Proposition
\ref{amplecritrelation}).

On the other hand, the existence of global syzygies also yields
results about the tight closure (Section \ref{primaryrelations}). If
there exists a primary syzygy (that is, a syzygy such that the
quotient is locally free) for $f_1,f_2,f_3$ of total degree $k \leq
(d_1+d_2+d_3)/2$, then $ f_0 \in (f_1,f_2,f_3)^\star$ holds for $m
\geq d_1+d_2+d_3-k$ and $f_0 \in (f_1,f_2,f_3)^\star $ if and only
if $f_0 \in (f_1,f_2,f_3)$ for $m <k$ (Corollary
\ref{primaryexclusion}).

In Section \ref{examplea} we consider the tight closure
$(x^a,y^a,z^a)^\star$ in $R=K[x,y,z]/(F)$, where $F$ is a
homogeneous polynomial of degree $\delta$ defining a smooth
projective curve. Our main result is that $(x^a,y^a,z^a)^\star =
(x^a,y^a,z^a) + R_{\geq \frac{3}{2}a} $ holds for $\delta \geq 3a-1$
(Corollary \ref{socloa}) in characteristic zero. In positive
characteristic we show that $ R_m \subset (x^a,y^a,z^a)^\star$ (and
even in the Frobenius closure) holds for $m > \frac{3}{2}a$, $
\delta \geq 3a-1$, $p > \delta -3$ and that $R_m \cap
(x^a,y^a,z^a)^\star = R_m \cap (x^a,y^a,z^a)$ holds for $m <
\frac{3}{2}a$, $ \delta \geq 3a-1$, $p \gg 0$ (Corollary
\ref{frobeniustight}).

In Section \ref{tightexample} we have a closer look at
$(x^2,y^2,z^2)^\star$ and $(x^3,y^3,z^3)^\star$ in $K[x,y,z]/(F)$
for $F$ of low degree $\delta$. We extend some of the general
results to degree $\delta < 3a-1$ and we prove some of the
conjectures of Sullivant to which he was led by his computations of
$(x^3,y^3,z^3)^\star$ in the Fermat rings.

I would like to thank the referee for careful reading and useful remarks.

\section{Slope criteria for tight closure}
\label{sectionrecall}

We recall the main results of \cite{brennerslope}. Let $R$ denote a
normal two-dimensional standard-graded $K$-domain over an
algebraically closed field $K$ and let $f_1, \ldots ,f_n$ denote
homogeneous $R_+$-primary elements of degree $d_i =\deg (f_i)$.
These elements define the locally free sheaf of syzygies $ \Syz(m)$
on the smooth projective curve $Y= \Proj \, R$ given by the short
exact sequence
$$0 \lra \Syz(m) \lra \bigoplus_{i=1}^n \O_Y(m-d_i) \stackrel{f_1, \ldots, f_n}
{\lra} \O_Y(m) \lra 0 \, .$$ Another homogeneous element $f_0$ of
degree $m$ yields via the connecting homomorphism a cohomology class
$c= \delta (f_0) \in H^1(Y, \Syz(m))$. This class corresponds to an
extension (set $\shS= \Syz(m)$) $c \in H^1(Y, \shS) \cong
Ext^1(\O_Y, \shS)$, say $0 \ra \shS \ra \shS' \ra \O_Y \ra 0$, or $0
\ra V \ra V' \ra \AA_Y \ra 0$ if written for geometric vector
bundles with sheaf of sections $\shS$. The main equivalence is now
that
$$(1.1)  \hspace{10pt} f_0 \in (f_1, \ldots ,f_n)^\soclo \hspace{10pt} \mbox{  if and only if  }
\hspace{10pt} \PP(V') - \PP(V) \hspace{6pt}  \mbox{ is not an affine
scheme} \, ,$$ where $\PP(V)$ denotes the projective bundle
corresponding to the vector bundle $V$ (this is $\PP(\shS^\dual)$ in
the notation of \cite{EGAII} or \cite[II.7]{haralg}). This
equivalence rests upon the interpretation of tight closure as solid
closure and the geometric interpretation of forcing algebras; see
\cite{hochstersolid} and \cite{brennertightproj}. We will work with
solid closure denoted by $I^\soclo$ in the following even if we
speak about tight closure. It is the same as tight closure in
positive characteristic.

For the affineness of $\PP(V') - \PP(V)$ given by a class $c \in
H^1(Y, \shS)$ we have proved in \cite{brennerslope} several
sufficient and necessary slope criteria. Recall that the slope of a
locally free sheaf $\shS$ is defined by $\mu (\shS)= \deg \, (\shS)
/ \rank (\shS)$. Every locally free sheaf $\shS$ has a unique {\em
Harder-Narasimhan filtration} (see \cite{hardernarasimhan},
\cite{huybrechtslehn}, \cite{lazarsfeldpositive}). This is a
filtration of locally free subsheaves
$$ 0=\shS_0 \subset \shS_1 \subset \ldots \subset \shS_s = \shS $$
such that $\shS_i/\shS_{i-1}$ is semistable for every $i=1, \ldots,
s$. Here $\shS_1$ is called the {\em maximal destabilizing
subsheaf}. The slopes of these semistable quotients form a
decreasing chain $\mu_1 > \ldots > \mu_s$. We call $\mu_{\rm
min}(\shS)= \mu_s= \mu(\shS/\shS_{s-1})$ the {\em minimal slope} and
$\mu_{\max}(\shS) =\mu_1 (\shS)$ the {\em maximal slope}. This is
the same as $\mu_{\rm min}(\shS)=\min\{\mu(\shQ):\, \shS \ra \shQ
\ra 0 \mbox{ locally free quotient sheaf} \}$ and $\mu_{\rm
max}(\shS) = \max \{\mu(\shT): \, \shT \subseteq \shS \mbox{ locally
free subsheaf} \}$. For the dual sheaf we have $\mu_{\rm max}
(\shS^\dual )= -\mu_{\rm min}(\shS)$. A locally free sheaf $\shS$ is
called semistable if $\mu_{\rm min}(\shS) =\mu_{\rm max}(\shS)$.

\medskip
In positive characteristic we need the following definition.

\begin{definition}
Let $Y$ denote a smooth projective curve over an algebrai\-cal\-ly closed
field and let
$\shS$ denote a locally free sheaf. Then we define
$$
\bar{\mu}_{\max} (\shS)=
{\rm sup} \{  \frac{\mu_{\rm max} (\varphi^* \shS)}{\deg \, (\varphi) }| \,\,
\varphi: Z \ra Y \mbox{ finite dominant $K$-morphism}      \}
$$
and
$$
\bar{\mu}_{\rm min} (\shS)=
{\rm inf} \{  \frac{\mu_{\rm min} (\varphi^* \shS)}{\deg \, (\varphi)} |\, \,
\varphi: Z \ra Y \mbox{ finite dominant $K$-morphism}  \} \, .
$$
\end{definition}

These numbers exist and give nothing new in characteristic zero. In
positive characteristic however they may differ from $\mu_{\rm
max}(\shS)$ and $\mu_{\rm min }(\shS)$. We say $\shS$ is {\em
strongly semistable} if $\bar{\mu}_{\rm max}(\shS) = \bar{\mu}_{\rm
min }(\shS)$. This is equivalent to the property that every
Frobenius pull-back of $\shS$ is semistable; see \cite[\S
5]{miyaokachern}.

With these notions our main slope criteria for affineness are the
following.

\begin{theorem}
\label{slopecriteria}
Let $Y$ denote a smooth projective curve over an algebraically closed field $K$,
let $\shS$ denote a locally free sheaf on $Y$ and let
$c \in H^1(Y, \shS)$ denote a cohomology class given rise
to $\PP(V') - \PP(V)$.
Then the following hold.

\renewcommand{\labelenumi}{(\roman{enumi})}
\begin{enumerate}

\item
Suppose that the characteristic of $K$ is zero.
If $\mu_{\rm max} (\shS)<0$ and $c \neq 0$, then $\PP(V') - \PP(V)$ is affine.

\item
Suppose that the characteristic of $K$ is zero.
Suppose that there exists a sheaf homomorphism
$\varphi: \shS \ra \shT$ such that $\shT$ is semistable of negative slope
and $0 \neq \varphi(c) \in H^1(Y, \shT)$.
Then $\PP(V') - \PP(V)$ is affine.

\item
If $\bar{\mu}_{\rm min}(\shS) \geq 0$, then $\PP(V') - \PP(V)$ is not affine.

\item
If $c=0$, then $\PP(V') - \PP(V)$ is not affine.

\end{enumerate}
\end{theorem}

\proof See \cite[Corollary 3.2, Theorem 3.4, and Theorem
4.4]{brennerslope} ((iv) is trivial). \qed

\medskip
The condition in Theorem \ref{slopecriteria}(i) implies that the
dual sheaf $\shS^\dual$ and also the extension $\shS'^\dual$ given
by $c \neq 0$ is ample (\cite[Theorem 2.2]{giesekerample}). This
means by definition that the divisor $\PP(V) \subset \PP(V')$ is
ample and hence its complement is affine. Ampleness and affineness
are open properties: if we have a smooth projective relative curve
$Y$ over $\Spec\, D$, where $\ZZ \subseteq D$ is a finitely
generated $\ZZ$-algebra, and if $\shS$ is locally free on $Y$, then
the affineness of $\PP(V'_\eta)- \PP(V_\eta)$ over the generic point
$\eta \in \Spec\, D$ implies the affineness of $\PP(V')- \PP(V)$
over an open non-empty subset of $\Spec\, D$. This observation
allows us to deduce from results in characteristic zero results for
characteristic $p \gg 0$, in particular when the situation is given
by a tight closure problem. From these affineness criteria we get
the following slope criteria for tight closure.

\begin{theorem}
\label{slopetightcriteria} Let $R$ denote a two-dimensional normal
standard-graded domain over an algebraically closed field $K$. Let
$(f_1, \ldots ,f_n)$ denote an $R_+$-primary homogeneous ideal given
by homogeneous ideal generators of degree $d_i= \deg (f_i)$. Let $Y
= \Proj\, R$ denote the corresponding smooth projective curve of
degree $\delta = \deg \O_Y(1)$. Let $\Syz(m)$ denote the locally
free sheaf of syzygies of total degree $m$. Set $\mu_{\rm max}(f_1,
\ldots, f_n)\! := \! \mu_{\rm max} (\Syz(0)^\dual)$ and define
$\bar{\mu}_{\rm max}(f_1, \ldots, f_n)$, $\mu_{\rm min}(f_1, \ldots,
f_n)$ and $ \bar{\mu}_{\rm min}(f_1, \ldots, f_n)$ in the same way.
Let $f_0$ denote another homogeneous element. Then the following
hold.

\renewcommand{\labelenumi}{(\roman{enumi})}
\begin{enumerate}

\item
If $\deg(f_0) \geq \bar{\mu}_{\rm max} (f_1, \ldots ,f_n)/ \delta $,
then $f_0 \in (f_1, \ldots ,f_n)^\soclo$.

\item
Suppose that the characteristic of $K$ is zero or $p \gg 0$.

\noindent
If $\deg (f_0) < {\mu}_{\rm min} (f_1, \ldots ,f_n)/ \delta $,
then $f_0 \in (f_1, \ldots ,f_n)^\soclo $
if and only if $f_0 \in (f_1, \ldots ,f_n)$.

\item
Suppose that the characteristic is zero. If the sheaf of syzygies
$\Syz(m)$ is semistable, then
$$ (f_1, \ldots ,f_n)^\soclo =
(f_1, \ldots ,f_n) + R_{\geq \frac{d_1 + \ldots +d_n}{n-1}} \, .$$
\end{enumerate}
\end{theorem}
\proof These statements follow from Theorem \ref{slopecriteria}; see
\cite[Theorem 6.4, Theorem 7.3 and Theorem 8.1]{brennerslope} \qed

\begin{remark}
Note that $\det \Syz(m)= \O_Y((n-1)m- d_1 - \ldots -d_n)$ and that
$\deg (\Syz(m)) =( (n-1)m- d_1 - \ldots -d_n) \delta$, where $\delta
$ is the degree of the curve $Y$. Therefore
$$ \mu_{\rm min}(f_1, \ldots ,f_n) \leq
\mu(\Syz(0)^\dual)= \frac{d_1+ \ldots +d_n}{n-1} \delta \leq
\mu_{\rm max}(f_1, \ldots ,f_n) \, .$$
\end{remark}

For the actual computation of tight closure we have to find bounds
for the minimal and the maximal slope for the sheaf of syzygies and
criteria for semistability.

\section{An algorithm for low rank}
\label{algorithmsection}

In this section we describe the first steps of an ``algorithm'' to
decide whether an open subset $\PP(V') -\PP(V)$ given by a
cohomology class $c \in H^1(Y, \shS)$ is affine or not, where $\shS$
is the sheaf of sections in the geometric vector bundle $V$. It
always gives a complete answer if the maximal destabilizing subsheaf
of $\shS$ is semistable, hence in particular if the rank of $\shS$
is two. This implies that it is possible to decide whether $f_0 \in
(f_1,f_2,f_3)^\soclo$ holds or not, at least if we are able to
compute the Harder-Narasimhan filtration of the sheaf of syzygies.
We assume that the characteristic of $K$ is zero. The
Harder-Narasimhan filtration of $\shS$,
$$0 = \shS_0 \subset \shS_1 \subset \ldots \subset
\shS_{s-1} \subset \shS_s =\shS \, ,$$
splits into short exact sequences
$$0 \lra \shS_{j-1} \stackrel{i}{\lra} \shS_j
\lra \shS_j/\shS_{j-1} \lra 0 \, ,$$ where the quotients $\shS_j
/\shS_{j-1}$ are semistable with slope $\mu_j (\shS)= \mu(\shS_j
/\shS_{j-1})$. The algorithm uses the fact that for a cohomology
class $c_j \in H^1(Y, \shS_j)$ we have either $0 \neq \bar{c}_j \in
H^1(Y, \shS_j/\shS_{j-1})$ or $c_j = i(c_{j-1})$, where $c_{j-1} \in
H^1(Y,\shS_{j-1})$. This argumentation scheme requires arbitrary
subsheaves, so even if we start with a sheaf of syzygies
$\shS=\Syz(f_1, \ldots ,f_n)(m)$ and a cohomology class $c \in
H^1(Y, \Syz(m))$ given by another homogeneous element, arbitrary
sheaves and cohomology classes come naturally into play.

If $s=1$, then $\shS$ is semistable and everything is clear
by Theorem \ref{slopecriteria}(i), (iii).
If $s=2$, then we have an exact sequence
$$0 \lra \shS_1 \lra \shS \lra \shQ \lra 0 \, ,$$
where $\shS_1$ and $\shQ$ are semistable of different slope.
In this case the algorithm gives a complete answer.

We present the algorithm in the following diagram. Note that for
$s=2$ we have $\mu_{s-1} (\shS)= \mu_1(\shS) = \mu_{\rm max}
(\shS)$. So if this number is $<0$, then we may conclude that
$\PP(V') - \PP(V)$ is affine by Theorem \ref{slopecriteria}(i).

\newpage

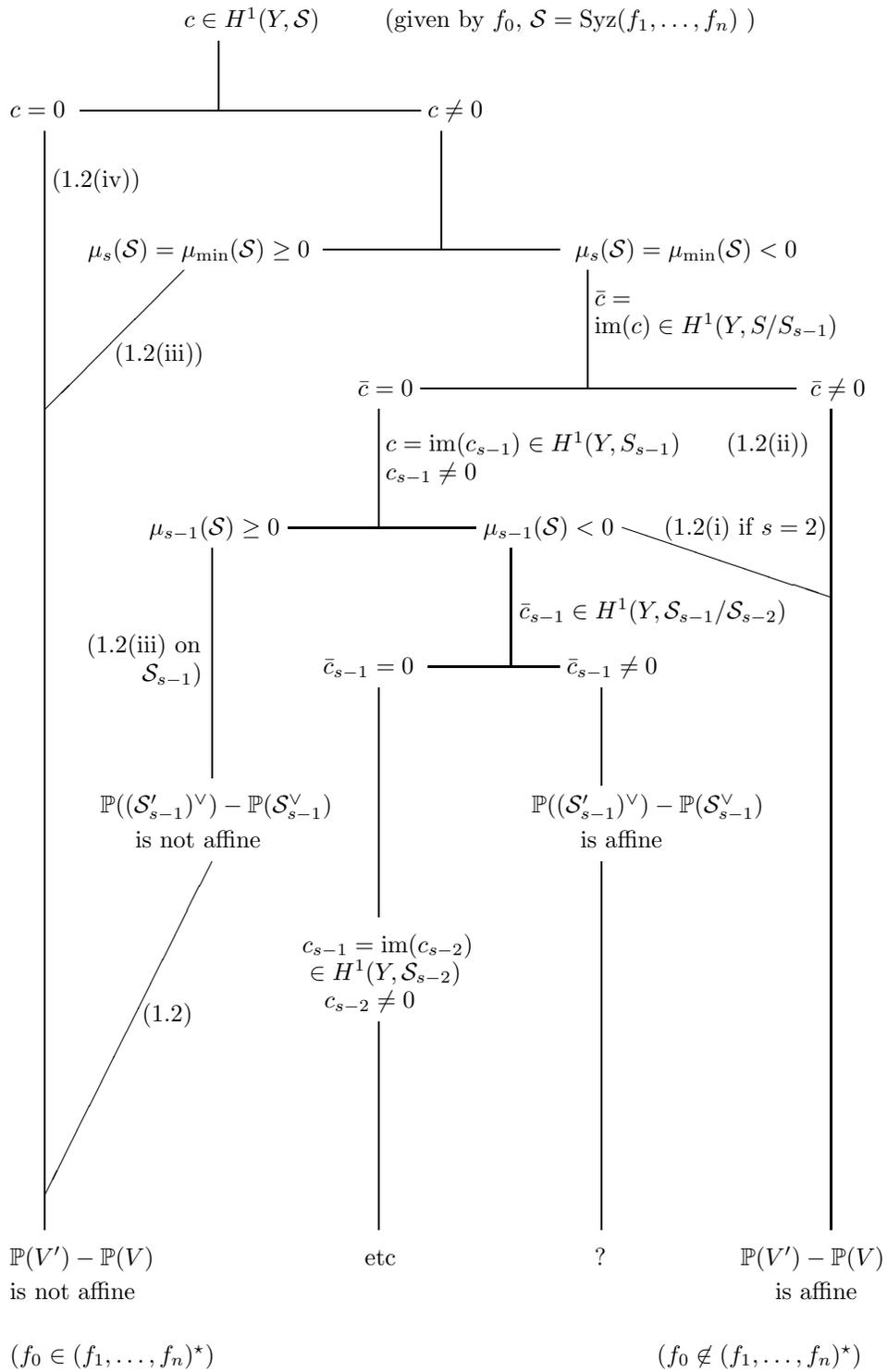
\small

\thispagestyle{empty}

\setlength{\unitlength}{1cm}

\begin{picture}(12,20.5)

\put(2.5, 20.3){$c \in H^1(Y, \shS) \mbox{ \ \ \ \ \ \ }$ (given by
$f_0$, $\shS= \Syz (f_1, \ldots ,f_n)$ )}

\put(3,20.1){\line(0,-1){1}}

\put(0,19){$c=0$}
\put(6,19){$c \neq 0$}
\put(1,19.1){\line(1,0){4.9}}

\put(0.5, 18.8){\line(0,-1){15.8}}

\put(0.6, 18){\rm (\ref{slopecriteria}(iv))}

\put(6.2 ,18.8){\line(0,-1){1.7}}

\put(1 ,17){ $\mu_s(\shS) = \mu_{\rm min}(\shS) \geq 0 $}

\put(8 ,17){ $\mu_s(\shS) = \mu_{\rm min}(\shS) < 0 $}

\put(4.5 ,17.1){\line(1,0){3.4}}

\put(2.5,16.8){\line(-1,-1){2.0}}

\put(8.3 ,16.8){\line(0,-1){1.7}}

\put(8.4 ,16.3){$ \bar{c}= $}

\put(8.4 , 15.9){$ {\rm im}(c) \in H^1(Y, S/ S_{s-1}) $}

\put(1.5, 15.5){\rm (\ref{slopecriteria}(iii))}

\put(5,15){$\bar{c} =0$}

\put(11.5 ,15){$ \bar{c} \neq 0$}

\put(5.9, 15.1){\line(1,0){5.4}}

\put(5.3 ,14.8){\line(0,-1){1.7}}

\put(11.8 ,14.8){\line(0,-1){11.8}}

\put(5.4 ,14.2){$ c = {\rm im}(c_{s-1}) \in H^1(Y, S_{s-1}) $}

\put(5.4, 13.8){$ c_{s-1} \neq 0$}

\put(10.3,14.2){\rm (\ref{slopecriteria}(ii))}

\put(2,13){$ \mu_{s-1}(\shS) \geq 0$}

\put(6.8,13){$ \mu_{s-1}(\shS) < 0$}

\put(4, 13.1){\line(1,0){2.7}}

\put(2.9,12.8){\line(0,-1){3.3}}

\put(1.1, 11.3){(\rm \ref{slopecriteria}(iii) on}

\put(1.9, 10.9){$\shS_{s-1}$)}

\put(7.2 ,12.8){\line(0,-1){1.7}}

\put(8.8 ,13.1){\line(3,-1){3}}

\put(9.4,13){\rm (\ref{slopecriteria}(i) if $s=2$)}

\put(7.2, 11.8){ $ \bar{c}_{s-1} \in H^1(Y,\shS_{s-1}/ \shS_{s-2})$}

\put(4.5,11){$\bar{c}_{s-1} =0$}

\put(8 ,11){$ \bar{c}_{s-1} \neq 0$}

\put(6. ,11.1){\line(1,0){1.9}}

\put( 8.5, 10.8){\line(0,-1){1.4}}

\put(5.3 ,10.8){\line(0,-1){3.3}}

\put(1.3, 9){$ \PP(( \shS_{s-1}')^\dual) - \PP( \shS_{s-1}^\dual)$}

\put(1.8,8.5){is not affine}

\put(7.5, 9){$ \PP(( \shS_{s-1}')^\dual) - \PP( \shS_{s-1}^\dual)$}

\put(8.2,8.5){is affine}

\put(8.5 ,8.3){\line(0,-1){5.3}}

\put(2.9,8.3){\line(-1,-2){2.4}}

\put(1.9,6){\rm (\ref{slopecriteria})}

\put(4.2 ,7){$c_{s-1}= \im(c_{s-2})$}

\put(4.3 ,6.6){$ \in H^1(Y,\shS_{s-2})$}

\put(4.4,6.2){  $c_{s-2} \neq 0$ }

\put(5.3, 6.0){\line(0,-1){3}}

\put(0,2.5) {$ \PP(V') - \PP(V)$}

\put(0,2){is not affine}

\put(5.1 ,2.5){etc}

\put(8.4 ,2.5){?}

\put(10.5,2.5){$ \PP(V') - \PP(V)$}

\put(11,2){is affine}

\put(0,1.1){($ f_0 \in (f_1, \ldots, f_n)^\soclo $)}
\put(9.3,1.1){($ f_0 \not\in (f_1, \ldots, f_n)^\soclo $)}
\end{picture}

\normalsize

\newpage

\section{The tight closure of three elements}
\label{tightthree}

We fix the following situation.

\begin{situation}
\label{situation} Let $R$ denote a two-dimensional normal
standard-graded domain over an algebraically closed field $K$. Let
$f_1, f_2, f_3 \in R$ denote three homogeneous $R_+$-primary
elements of degree $d_i = \deg (f_i)$. Let $Y= \Proj \, R$ denote
the corresponding smooth projective curve of degree $\delta = \deg
\O_Y(1)$ and of genus $g$ and let $\Syz(m)$ denote the sheaf of
syzygies of total degree $m$ for $f_1,f_2,f_3$.
\end{situation}

The sheaf of syzygies $\Syz(m)$ on $Y$ has rank two; hence we may
decide in characteristic zero due to Section \ref{algorithmsection}
whether $\PP(V') - \PP(V)$ (given by the cohomology class
$\delta(f_0)=c \in H^1(Y, \Syz(m))$) is affine or not. Therefore we
may compute $(f_1,f_2,f_3)^\soclo$, at least if we can compute the
Harder-Narasimhan filtration of $\Syz$.

In positive characteristic we have to refine this algorithm,
since the Harder-Nara\-simhan filtration is in general not stable
under the Frobenius pull-back. For a locally free sheaf $\shS$ we set
$$ \maxdeg_1 (\shS)
= \max \{ \deg(\shL):\, \shL
\mbox{ is an invertible subsheaf of } \shS \} $$
and
$$ \mindeg_1 (\shS)
= \min \{ \deg(\shL):\, \shL \mbox{ is an invertible quotient sheaf
of } \shS \} \, .$$ For a locally free sheaf $\shS$ of rank two we
have $\mu_{\rm max} (\shS) = \max \{ \maxdeg_1 (\shS) ,\deg (\shS)/2
\}$ and $\mu_{\rm min} (\shS) = \min \{ \mindeg_1 (\shS), \deg
(\shS)/2 \} $, and $\mindeg_1 (\shS)= \deg (\shS) -
\maxdeg_1(\shS)$. The sheaf $\shS$ is semistable if and only if
$\maxdeg_1 (\shS) \leq \deg (\shS)/2$.

If we find a subsheaf $\shL \subseteq \shS$ such that $\deg (\shL)
\geq \deg (\shS)/2$ and such that the quotient is itself locally
free (i.e., $\shL$ is a subbundle), then $\deg (\shL) =
\maxdeg_1(\shS) = \mu_{\rm max} (\shS)$.

\begin{lemma}
\label{exactsequenceranktwo}
Let $\shS$ denote a locally free sheaf of rank two on
a smooth projective curve $Y$.
Suppose that we have a short exact sequence
$0 \ra \shL \ra \shS \ra \shM \ra 0$, where $\shL$ and $\shM$ are invertible
sheaves.
Then $\maxdeg_1 (\shS) \leq \max (\deg (\shL), \deg (\shM) )$.
If furthermore $\deg (\shL)  \geq \mu (\shS) = \deg (\shS)/2$,
then $\maxdeg_1 (\shS) = \deg (\shL)$
and $\mindeg_1(\shS) = \deg (\shM)$.

If $\deg (\shL) = \mu(\shS)$, then $\shS$ is strongly semistable.
If $\deg (\shL) > \mu(\shS)$, then $\shS$ is not stable
and $\shL$ is the maximal destabilizing subsheaf.
\end{lemma}

\proof Let $\shN$ denote an invertible sheaf and let $\varphi: \shN
\ra \shS$ be a sheaf morphism. If $\deg (\shN) > \deg (\shM)$, then
the composed morphism $ \shN \ra \shM$ is zero and $\varphi$ factors
through $\shL$. But then $\deg (\shN) \leq \deg (\shL)$ or $\varphi$
is zero. So suppose that $\deg (\shL)  \geq \mu (\shS)$. Then $\deg
(\shM) = \deg (\shS) -\deg (\shL) \leq \deg (\shS) - \deg (\shS)/2 =
\deg (\shS)/2$, and hence $ \deg (\shL) \geq \mu ( \shS) \geq \deg
(\shM)$. Thus $\shL$ is a subbundle of maximal degree and $\shM$ is
a quotient invertible sheaf of minimal degree. The other statements
follow. \qed

\medskip
We now have the following two alternatives. The locally free sheaf
$\shS$ of rank two on $Y$ is strongly semistable. Then
$\bar{\mu}_{\rm max}(\shS) = \mu(\shS)= \bar{\mu}_{\rm min}(\shS)$
and the tight closure is easy to compute by the numerical criterion
Theorem \ref{slopetightcriteria}. Otherwise $\shS$ is not strongly
semistable. Then there exists a finite morphism $\varphi: Y' \ra Y$
such that there exists a short exact sequence on $Y'$, $ 0 \ra \shL
\ra \shS' \ra \shM \ra 0$, where $\deg (\shL) \geq \mu (\shS')$,
$\shS' = \varphi^*(\shS)$. In this case the pull-back of this
sequence for another morphism $\psi : Y'' \ra Y'$ also fulfills the
condition in Lemma \ref{exactsequenceranktwo}; hence $\mu_{\rm max}
(\shS'') = \deg (\psi^*(\shL))= \deg (\shL) \deg (\psi)$ and thus
$\bar{\mu}_{\rm max} (\shS) = \deg (\shL)/ \deg (\varphi)$.

\medskip
If we have a short exact sequence for the sheaf of syzygies for
three elements $f_1,f_2,f_3$, then we can often compute
$(f_1,f_2,f_3)^\soclo$ and $(f_1,f_2,f_3)^+$ (the plus closure)
according to the following proposition.

\begin{proposition}
\label{exactsequencecrit} Suppose the notation and situation of {\rm
\ref{situation}}. Let $0 \ra \shL \ra \Syz(m) \ra \shM \ra 0$ denote
a short exact sequence, where $\shL$ and $\shM$ are invertible. Let
$f_0$ denote another homogeneous element of degree $m$ and let $c
\in H^1(Y, \Syz(m))$ denote its forcing class. Let $\bar{c}$ denote
the image of $c$ in $H^1(Y, \shM)$. Then the following hold {\rm
(}suppose in the first two statements that the characteristic of $K$
is zero or $p \gg 0${\rm )}.

\renewcommand{\labelenumi}{(\roman{enumi})}
\begin{enumerate}

\item
If $\deg (\shL) <0 $ and $\deg (\shM) <0 $ and $c \neq 0$,
then $f_0 \not\in (f_1,f_2,f_3)^\soclo$.

\item
If $\deg (\shM) <0 $ and
$\bar{c} \neq 0$, then $f_0 \not\in (f_1,f_2,f_3)^\soclo$.

\item
If $\bar{c} =0$ and $\deg (\shL) \geq 0$, then
$f_0 \in (f_1,f_2,f_3)^\soclo$.

\item
If $\deg (\shL) \geq \deg (\shM) \geq 0 $,
then $f_0 \in (f_1,f_2,f_3)^\soclo$.

\item
Suppose that the characteristic of $K$ is positive.
Suppose that $\bar{c}=0 $ or that $\shM$ is trivial or that $\shM$ has
positive degree. Moreover suppose that $\shL$ is trivial or has positive degree.
Then $f_0 \in (f_1,f_2,f_3)^+$.

\end{enumerate}
\end{proposition}

\proof We use the equivalence that $f_0 \not\in
(f_1,f_2,f_3)^\soclo$ if and only if $\PP(V') - \PP(V)$ is affine
(1.1), where $V$ is the vector bundle corresponding to the sheaf of
syzygies and $V'$ is the extension of $V$ given by the forcing
class. (i) The maximal slope of $\Syz(m)$ is $< 0$; hence the result
follows from Theorem \ref{slopecriteria}(i).

(ii) This follows from \ref{slopecriteria}(ii).

(iii) If $\bar{c}=0$, then there exists a cohomology class $e \in
H^1(Y, \shL)$ mapping to $c$. Let $\shL'$ denote the extension of
$\shL$ defined by this cohomology class $e$. If $f_0 \not \in
(f_1,f_2,f_3)^\soclo$, then $\PP(V') - \PP(V)$ would be affine and
then $\PP((\shL') ^\dual) - \PP(\shL^\dual) \subset \PP(V') -
\PP(V)$ would be affine as a closed subscheme, but this is not true
since $\deg (\shL) \geq 0$ and Theorem \ref{slopecriteria}(iii).

(iv) The condition implies in connection with Theorem
\ref{exactsequenceranktwo} that $\bar{\mu}_{\rm min}( \Syz(m)) \geq
0$; hence the result follows again from Theorem
\ref{slopecriteria}(iii).

(v) After applying a finite mapping $\varphi: Y' \ra Y$, we may
assume that $\bar{c} =0$. For $\deg (\shM) >0$ this can be done by a
Frobenius power and for $\shM = \O_Y$ this is due to
\cite[Proposition 8.1]{brennertightproj}. Therefore we may assume
that $\varphi^*(c)$ stems from a cohomology class $e \in H^1(Y',
\varphi^*(\shL))$. Due to the assumptions on $\shL$ we can do the
same with $e$; hence there exists altogether a finite mapping $Y''
\ra Y$ such that the pull-back of $c$ is zero. Therefore $f_0 \in
(f_1,f_2,f_3)^+$. \qed

\medskip
We may apply Proposition \ref{exactsequencecrit} to the short exact
sequence given by the Harder-Narasim\-han filtration to compute the
tight closure of $(f_1, f_2, f_3)$ if $\Syz$ is not strongly
sta\-ble.

\begin{corollary}
\label{exactsequencedecide} Suppose the situation and notation of
{\rm \ref{situation}}. Suppose that the sheaf of syzygies $\Syz$ is
not strongly stable, and let $\varphi: Y' \ra Y$ denote a finite
dominant morphism of smooth projective curves such that there exists
a short exact sequence
$$0 \lra \shL(m) \lra \varphi^*(\Syz(m)) \lra \shM(m) \lra 0$$ on $Y'$,
where $\deg (\shL(m)) \geq \mu (\varphi^*(\Syz (m)) ) \geq \deg
(\shM(m))$. Let $f_0$ be another homogeneous element of degree $m$,
let $c$ denote its forcing class in $H^1(Y',\varphi^*(\Syz(m)) )$
and let $\bar{c}$ denote its image in $H^1(Y, \shM(m))$. Then we may
decide whether $f_0 \in (f_1,f_2,f_3)^\soclo$ in the following way
{\rm(}assume in the first and second statement that the
characteristic is zero or $p \gg 0${\rm )}.

\renewcommand{\labelenumi}{(\roman{enumi})}
\begin{enumerate}

\item
If $\deg (\shL(m)) <0 $ and $c \neq 0$, then $f_0 \not\in (f_1,f_2,f_3)^\soclo$.

\item
If $\deg (\shL(m))  \geq 0 $, $\deg (\shM(m)) <0$
and $\bar{c} \neq 0$, then $f_0 \not\in (f_1,f_2,f_3)^\soclo$.

\item
If $\deg (\shL(m))  \geq 0 $, $\deg (\shM(m)) <0$
and $\bar{c} = 0$, then $f_0 \in (f_1,f_2,f_3)^\soclo$.

\item
If $\deg (\shM(m)) \geq 0 $, then $f_0 \in (f_1,f_2,f_3)^\soclo$.

\end{enumerate}
\end{corollary}

\proof This follows from Proposition \ref{exactsequencecrit}. \qed

\begin{remark}
\label{exactsequencepositive} Suppose the situation of
\ref{exactsequencedecide} and suppose that the characteristic is
positive. Then we need in (iii) and (iv) stronger conditions to
conclude that $f_0 \in (f_1,f_2,f_3)^{+}$. For (iii) we need that
$\shL(m)$ is trivial or that $\deg (\shL(m)) $ is positive. For (iv)
we need that $\shM$ is of positive degree. This follows from
Proposition \ref{exactsequencecrit}(v).
\end{remark}

If the sheaf of syzygies is decomposable, that is, the sum of two
invertible sheaves, then the decomposition gives at once a short
exact sequence and Proposition \ref{exactsequencecrit} and Corollary
\ref{exactsequencedecide} are particularly easy to apply. On the
other hand we cannot expect any bound for $\mu_{\rm max}$ in the
decomposable case. For the indecomposable case we have the following
result.

\begin{theorem}
\label{genusbound} Suppose the notation and situation of {\rm
\ref{situation}} and suppose that the characteristic of $K$ is zero.
Suppose that the sheaf of syzygies $\Syz(m)$ is indecomposable on
$Y$. Then
$$\mu_{\rm max} (f_1,f_2,f_3) \leq \delta \frac{d_1+d_2+d_3}{2} + g-1 \, $$
and
$$ \mu_{\rm min} (f_1,f_2,f_3) \geq \delta \frac{d_1+d_2+d_3}{2} - g+1 \, .$$
\end{theorem}

\proof
See \cite[Theorem 9.1]{brennerslope}.
\qed

\section{The degree of syzygies}
\label{tightrelations}

The notions of semistability and of minimal and maximal degree of a
locally free sheaf $\shS$ on a smooth projective curve $Y$ refer to
all locally free subsheaves of $\shS$ (or quotient sheaves). However
for a syzygy sheaf $\Syz(m)$ defined by homogeneous primary elements
$f_1,\ldots ,f_n$ in a two-dimensional normal standard-graded
$K$-domain $R$ we have the fixed polarization $\O_Y(1)$ on $Y=
\Proj\, R$. It is then often easier to control the behavior of
$\Syz(m) = \Syz(0) \otimes \O_Y(m)$ instead of $\Syz(0) \otimes
\shL$ for all invertible sheaves $\shL$. The (non-)\-existence of
syzygies $\neq 0$ for $f_1, \ldots ,f_n$ of certain degree has many
consequences on the structure of $\Syz(m)$ and hence on the
corresponding tight closure problem.

\begin{remark}
\label{riemannroch} The global syzygies of $\Syz(m)$ obey the
theorem of Riemann-Roch, that is,
$$ h^0(Y, \Syz(m)) - h^1(Y, \Syz(m))
=  ((n-1)m-d_1- \ldots -d_n) \delta+ (1-g) (n-1) \, ,$$ where
$\delta = \deg \O_Y(1)$ is the degree of $Y$ and $g$ is its genus.
In particular, for $ m > \frac{d_1+ \ldots + d_n}{n-1} +
\frac{g-1}{\delta}$ there exist global syzygies $\neq 0$ in
$\Syz(m)$.
\end{remark}

\begin{lemma}
\label{maxdeglemma}
Suppose that the locally free sheaf $\shS$ on the smooth
projective curve $Y$ of genus $g$ over an algebraically closed field $K$
does not have sections $\neq 0$.
Then $\Gamma(Y,\shS \otimes \shL) =0 $
for every invertible sheaf $\shL$ of degree $\leq -g$.
In particular $\maxdeg_1 (\shS) \leq g-1$.
\end{lemma}
\proof Suppose the contrary. Then we have a non-trivial morphism
$\shM \ra \shS$ such that $\deg \shM \geq g$. But due to the theorem
of Riemann-Roch we have $ h^0(\shM) \geq \deg (\shM) +1-g$; hence
the invertible sheaf $\shM$ must have non-trivial sections, which
gives a contradiction. \qed

\begin{proposition}
\label{maxinvertbound}
Let $f_1, \ldots, f_n$ denote homogeneous primary elements
in a normal two-dimensional standard-graded $K$-domain $R$ of degree $d_i$,
where $K$ is an algebraically closed field.
Suppose that $Y= \Proj\, R$ has genus $g=g(Y)$ and degree $\delta$.

\renewcommand{\labelenumi}{(\roman{enumi})}

\begin{enumerate}

\item
Suppose that there exists a syzygy $\neq 0$ for the elements $f_1,
\ldots, f_n$ of total degree $k < (d_1 + \ldots +d_n )/(n-1)$. Then
the sheaf of syzygies is not semistable.

\item
Suppose that there does not exist a syzygy $\neq 0$ of total degree
$k$. Then
$$ \maxdeg_1(\Syz(k)) \leq g-1 \, . $$

\item
Let $n=3$. Suppose that there exists a syzygy $\neq 0$ of total
degree $k < \frac{d_1 + d_2 +d_3}{2} - \frac{g-1}{\delta} $. Then
the sheaf of syzygies is decomposable, i.e., the sum of two
invertible sheaves.

\item
Let $n=3$ and suppose that there does not exist a syzygy $\neq 0$ of
total degree $k \geq \frac{d_1+d_2+d_3}{2} + \frac{g-1}{ \delta} $.
Then $\Syz$ is semistable.
\end{enumerate}
\end{proposition}

\proof (i) The syzygy $\neq 0$ induces a non-trivial morphism $\O_Y
\ra \Syz(k)$, but the degree $\deg (\Syz(k))= ((n-1)k - d_1- \ldots
-d_n) \delta < 0$ is negative; hence $\Syz(k)$ is not semistable.

(ii) The assumption means that $\shS= \Syz(k)$ has no global
sections $\neq 0$; hence Lemma \ref{maxdeglemma} yields that
$\maxdeg_1 (\Syz(k)) \leq g-1$.

(iii) Since $\Syz(k)$ has a non-trivial section, the sheaf $\Syz(k)$
contains the structure sheaf as a subsheaf and therefore $\mu_{\rm
max} (\Syz(k)) \geq 0$. On the other hand we have $\mu (\Syz(k))
+g-1 = (k- \frac{d_1+d_2+d_3}{2}) \delta +g-1 < 0$, that is,
$\mu(\Syz(k))<1-g$. Therefore we have a short exact sequence $0 \ra
\shL \ra \Syz(k) \ra \shM \ra 0$ where $\deg (\shL) \geq 0$ and
$\deg (\shM) < 2(1-g)$. This extension corresponds to a class
$H^1(Y, \shL \otimes \shM^{-1}) \cong H^0(Y, \shL^{-1} \otimes \shM
\otimes \omega_Y ) =0$; hence it is trivial.

(iv) The numerical condition means that $g-1 \leq   (k -
\frac{1}{2}(d_1+d_2+d_3)) \delta = \deg (\Syz(k))/2$; hence from
(ii) we get that $\maxdeg_1 (\Syz(k)) \leq \deg (\Syz(k))/2 $ and
the sheaf of syzygies is semistable. \qed

\medskip
We may derive from Proposition \ref{maxinvertbound}(ii)
the following inclusion bound for tight closure.

\begin{corollary}
\label{relationboundinclusion} Suppose the notation and situation of
{\rm \ref{situation}} and that the characteristic of $K$ is zero.
Suppose that there does not exist a syzygy $\neq 0$ of total degree
$k \leq \frac{d_1+d_2+d_3}{2} + \frac{g-1}{\delta}$. Then $R_m
\subseteq (f_1,f_2,f_3)^\soclo$ holds for $m \geq  d_1+d_2+d_3-k +
\frac{g-1}{\delta }$.
\end{corollary}
\proof If $\Syz$ is semistable, then the result follows from Theorem
\ref{slopetightcriteria}(iii). If $\Syz$ is not semi\-sta\-ble, then
$ \mu_{\rm max} (\Syz(k)) = \maxdeg_1 (\Syz(k)) \leq g-1 $ by
Proposition \ref{maxinvertbound}(ii) and we get
\begin{eqnarray*}
\mu_{\rm max} (f_1, f_2 ,f_3) &=& (d_1+d_2+d_3) \delta - \mu_{\rm
min} (f_1,f_2,f_3) \cr &=& (d_1+d_2+d_3) \delta + \mu_{\rm max}(
\Syz(0)) \cr &=& (d_1+d_2+d_3) \delta + \mu_{\rm max}( \Syz(k)) -k
\delta \cr &\leq & (d_1 +d_2+d_3 -k) \delta + g-1 \, .
\end{eqnarray*}
Hence under the numerical condition we have $\deg (f_0)=m \geq
\mu_{\rm max}(f_1, f_2,f_3)/ \delta $ and the result follows from
Theorem \ref{slopetightcriteria}(i). \qed

\begin{corollary}
\label{semistablecritcor} Let $R=K[x,y,z]/(F)$ denote a normal
two-dimensional stan\-dard-graded $K$-domain over an algebraically
closed field $K$, where $F$ is a polynomial of degree $\delta$. Let
$f_1,f_2,f_3 \in R$ be $R_+$-primary homogeneous elements of degree
$d_1,d_2,d_3$. Suppose that there does not exist a syzygy $\neq 0$
for $f_1,f_2,f_3$ of total degree $k$ with $k \geq
\frac{d_1+d_2+d_3}{2} + \frac{\delta -3}{2}$. Then $\Syz$ is
semistable.
\end{corollary}
\proof This follows directly from Proposition
\ref{maxinvertbound}(iv). \qed

\begin{remark}
Corollary \ref{semistablecritcor} is only applicable for $k=
\frac{d_1+d_2+d_3}{2} + \frac{\delta -3}{2}$, since for greater $k$
there exist global syzygies due to Remark \ref{riemannroch}. The
same is true for Proposition \ref{maxinvertbound}(iv).
\end{remark}

\begin{example}
\label{powersexample} We consider the elements $x^d,y^d,z^d$ on a
smooth projective curve given by an equation $F=0$, where $F$ is a
homogeneous polynomial of degree $\delta$. There exist syzygies like
$(y^d,-x^d,0)$ of total degree $2d$. Suppose that there do not exist
syzygies of smaller degree. Then the numerical condition in
Corollary \ref{semistablecritcor} for semistability is that $2d-1
\geq 3d/2 + (\delta-3)/2$ or equivalently that $\delta \leq d+1$. If
we want to apply Corollary \ref{semistablecritcor} we have to make
sure that the defining polynomial $F$ of degree $\delta$ does not
yield syzygies of degree $< 2d$.

\smallskip
Look at $d=2$ and $\delta =3$. If the monomial $xyz$ does occur in
$F$, then there do not exist syzygies of degree $3$ and the syzygy
sheaf is semistable. However this yields nothing interesting for
tight closure, since then $xyz \in (x^2,y^2,z^2)$ holds anyway.

\smallskip
Now look at $d=4$ and $\delta =5$. Under suitable conditions for the
coefficients of $F$ there does not exist a syzygy of degree $7$ for
$x^4,y^4,z^4$. Write $F=a x^3y^2+ bx^3yz+cx^3z^2+ dx^2y^3 + \ldots
$. A syzygy of degree $7$ is the same as a multiple $FQ$ ($\deg
(Q)=2$) which belongs to $(x^4,y^4,z^4)$. The six monomials of
degree $2$ yield six linear combinations in the six monomials of
degree $7$ outside $(x^4,y^4,z^4)$, namely $ x^3y^3z$, $x^3y^2z^2$,
$x^3yz^3$, $x^2y^3z^2$, $x^2y^2z^3$ and $xy^3z^3$. We may choose the
coefficients of $F$ in such a way that these linear combinations are
linearly independent. Then there does not exist a syzygy of degree
$7$. So in this case the sheaf of syzygies is semistable. It follows
for $\Char (K)=0$ that $R_6 \subseteq (x^4,y^4,z^4)^\soclo$. Note
that it is not true that $R_6 \subseteq (x^4,y^4,z^4)$, since there
exist $10$ monomials of degree $6$ outside $(x^4,y^4,z^4)$ in
$K[x,y,z]$; hence the dimension of $R_6/ (x^4,y^4,z^4)$ is at least
$10-3$.
\end{example}

\begin{remark}
\label{dual} If $\shS$ is a locally free sheaf of rank two, then we
have the natural mapping $\shS \oplus \shS \ra \shS \wedge \shS
\cong \det \shS$. This mapping induces an isomorphism $\shS \cong
\shS^\dual \otimes \det \shS$. If $\shS= \Syz(m)$ is the sheaf of
syzygies for three homogeneous elements $f_1,f_2,f_3$ of degree
$d_i$, then
\begin{eqnarray*}
(\Syz(m))^\dual \cong \Syz(m) \otimes (\det \Syz(m))^\dual & = &
\Syz(m) \otimes \O_Y( d_1+d_2+d_3-2m) \cr & =& \Syz(d_1+d_2+d_3-m)
\, . \end{eqnarray*} The natural mapping $\Syz(m) \oplus \Syz(m) \ra
\det (\Syz(m)) \cong \O_Y( 2m -d_1-d_2-d_3)$ is in terms of the
injection $\Syz(m) \subset \O_Y(m-d_1) \oplus \O_Y(m-d_2) \oplus
\O_Y(m-d_3)$ given by
$$(g_1,g_2,g_3), (h_1,h_2,h_3) \longmapsto
 \frac{ g_2h_3- g_3h_2}{f_1} = \frac{ -g_1h_3+g_3h_1}{f_2}
= \frac{ g_1h_2- g_2h_1}{f_3} \, .$$ A primary syzygy
$(g_1,g_2,g_3)$ $\in \Syz(m)$ yields the quotient mapping $\Syz(m)
\ra \O_Y(2m-d_1-d_2-d_3)$ which is given by $(h_1,h_2,h_3) \mapsto
\frac{-h_2g_3 +h_3g_2}{f_1}$.

If moreover $f_1$ and $f_2$ are parameters,
then the forcing class $c = \delta (f)$ of an element $f$ of degree $m$
is represented by the $\check{\rm C}$ech-cocycle
$$( \frac{f}{f_1}, - \frac{f}{f_2},0) \in H^1(Y, \Syz(m)) \, .$$
The quotient mapping sends this class to $- \frac{fg_3}{f_1f_2} \in H^1(Y, \O_Y(2m -d_1-d_2-d_3))$.
\end{remark}

\begin{example}
The existence of syzygies depends heavily on the characteristic.
Consider $x^{100}, y^{100}, z^{100}$ on $x^4+y^4+z^4=0$. For $\Char
(K) =5$ we find
$$(x^4+y^4+z^4)^{25} = x^{100} + y^{100} + z^{100} \, ;$$
hence there exists a syzygy of total degree $100$ and in fact
$z^{100} \in (x^{100}, y^{100})$ and therefore $(x^{100}, y^{100},
z^{100} )^\soclo=(x^{100}, y^{100})^\soclo$.

For $\Char (K)=37 $ we have $(x^4+y^4+z^4)^{37} = x^{48}x^{100} +
y^{48}y^{100} + z^{48}z^{100} $, which is a syzygy of total degree
$148$. This syzygy yields a short exact sequence
$$0 \lra \O_Y(4) \lra \Syz(152) \lra \O_Y \lra 0 \, ,$$
where the last mapping is given by $(g_1,g_2,g_3) \mapsto
\frac{x^{48}g_2- y^{48}g_1}{z^{100}}$ due to Remark \ref{dual}. Due
to Proposition \ref{maxinvertbound}(iii) the sheaf of syzygies is
decomposable, and we show that this sequence splits.

We may write $z^{100} = (-x^4-y^4)^{25} = x^{48}g_2 -y^{48}g_1$
where $g_1$ and $g_2$ have degree $52$. We can build a syzygy
$(g_1,g_2,g_3)$ of degree $152$ by $g_3 = \frac{-
(g_1x^{100}+g_2y^{100})}{z^{100}}$, for then
\begin{eqnarray*}
g_3  =  \frac{- (g_1x^{100}+g_2y^{100})}{z^{100}}
&=& \frac{- (g_1x^{148}+g_2y^{100}x^{48})}{x^{48}z^{100}} \cr
&=& \frac{ g_1(y^{148}+z^{148}) -g_2y^{100}x^{48}}{x^{48}z^{100}} \cr
&=& \frac{g_1 z^{48}}{x^{48}}
+ \frac{g_1y^{148}- g_2y^{100}x^{48}}{x^{48}z^{100}} \cr
&=&  \frac{g_1 z^{48}}{x^{48}} - \frac{y^{100} }{x^{48}}
\, .
\end{eqnarray*}
Therefore $g_3$ may be written with denominators $z$ and $x$; hence
it is a global section of $\O_Y(52)$. Hence $(g_1,g_2,g_3)$ is a
global syzygy of total degree $152$ which maps to $1$, so the short
exact sequence splits.
\end{example}

The nonexistence of global syzygies implies the ampleness of the
dual sheaf of the sheaf of syzygies. This observation then yields
exclusion criteria for tight closure.

\begin{proposition}
\label{amplecritrelation} Let $K$ denote an algebraically closed
field of characteristic zero and let $R$ denote a two-dimensional
normal standard-graded domain. Suppose that $f_1, \ldots , f_n \in
R$ are $R_+$-primary homogeneous elements. Suppose that there does
not exist a global syzygy $\neq 0$ of total degree $k$ for the
elements $f_1, \ldots , f_n$. Set $\shF(-m) := \Syz(m)^\dual$. Then
$\shF(-m)$ is ample for $m < k- \frac{2n-3}{(n-1) \delta} g +
\frac{1}{ \delta}$. In particular, $\shF(-m)$ is ample for $m \leq k
- 2g/\delta $.
\end{proposition}

\proof We know by Proposition \ref{maxinvertbound}(ii) that
$\maxdeg_1 (\Syz(k)) \leq g-1$ and therefore dually that $\mindeg_1
(\shF(-k)) \geq -g+1$. Hence
$$\mindeg_1 (\shF(-m)) = \mindeg_1 (\shF(-k) \otimes \O_Y(k-m))
= \mindeg_1 (\shF(-k)) + (k-m)\delta \geq -g+1  + (k-m) \delta \, .
$$ The numerical condition is equivalent to $ -g+1  + (k-m) \delta
> \frac{n-2}{n-1} g$. Hence $\mindeg_1(\shF(-m)) > \frac{n-2}{n-1}
g$ and the result follows from \cite[Corollary 2.2]{brennerslope}.
\qed

\begin{corollary}
\label{globalrelationample} Let $R=K[x,y,z]/(F)$ be a normal
standard-graded $K$-domain over an algebraically closed field $K$ of
characteristic zero, where $F$ is an irreducible polynomial of
degree $\delta$. Let $f_1, \ldots ,f_n$ denote primary homogeneous
elements of degree $d_i$. Suppose that there does not exist a global
syzygy $\neq 0$ for $f_1, \ldots, f_n$ of total degree $k$. Then
$\shF(-m)$ is ample for $m \leq k- \delta + 2$. An element $f_0 \in
R$ of degree $m \leq k- \delta + 2$ belongs to $(f_1, \ldots
,f_n)^\soclo$ if and only if it belongs to $(f_1, \ldots ,f_n)$.
\end{corollary}
\proof We have $m \leq k- \delta +2 \leq k - (\delta -2)(\delta
-1)/\delta = k - 2g/\delta $; hence the ampleness of $\shF(-m)$ for
$m \leq k - \delta +2$ follows from Proposition
\ref{amplecritrelation}. Now if $f$ is a homogeneous element of
degree $m \leq k -\delta +2$ such that $f \not\in (f_1 , \ldots ,
f_n)$, then $f$ defines a non-trivial extension $0 \ra \O \ra \shF '
\ra \shF(-m) \ra 0$ and therefore $\shF'$ is also ample
(characteristic zero). Therefore $\PP(\shF') - \PP(\shF(-m))$ is
affine and $f \not \in (f_1, \ldots, F_n)^\soclo $ by (1.1). \qed

\begin{example}
\label{powersexample2} We want to apply Proposition
\ref{amplecritrelation} and Corollary \ref{globalrelationample} to
Example \ref{powersexample} for $d=4$, $ \delta =5$ under the
condition that there does not exist a global syzygy of degree $7$.

Then Corollary \ref{globalrelationample} shows ($k=7$, $g=6$, $n=3$)
that $\shF(-m)$ is ample only for $m \leq 7-5+2=4$. The second bound
in Proposition \ref{amplecritrelation} gives this for $m \leq 7-
\frac{2 \cdot 6}{5} =4.6$. The first bound in Proposition
\ref{amplecritrelation} however yields ampleness for $m < 7 -
\frac{3}{2 \cdot 5}6 + \frac{1}{5} =7- \frac{16}{10}= 5.4 $.
Therefore $(x^4,y^4,z^4)^\soclo \cap R_m =(x^4,y^4,z^4) \cap R_m$
for $m \leq 5$.
\end{example}

\section{The existence of primary syzygies}
\label{primaryrelations}

We suppose further that $\Syz$ is the sheaf of syzygies on a smooth
projective curve $Y= \Proj R$ for homogeneous primary elements $f_1,
\ldots ,f_n \in R$, where $R$ is a two-dimensional normal
standard-graded $K$-domain over an algebraically closed field $K$.
We say that a syzygy $r \in \Gamma(Y,\Syz(m))$ is a {\em primary
relation} if it has no zero on $Y$ or, equivalently, if $r : \O_Y
\ra \Syz(m)$ defines a subbundle. For a primary syzygy we get a
short exact sequence $0 \ra \O_Y \ra \Syz(m) \ra \shQ \ra 0$, where
$\shQ$ is also locally free.

\begin{corollary}
\label{primarysemistable} Suppose the situation and notation of {\rm
\ref{situation}}. Suppose that there exists a primary syzygy of
total degree $k$. Then this syzygy gives rise to a short exact
sequence
$$0 \ra \O_Y \ra \Syz(k) \ra \O_Y(2k - d_1-d_2-d_3) \ra 0 \, .$$
If $k \geq \frac{d_1+d_2+d_3}{2}$,
then $R_m \subseteq (f_1,f_2,f_3)^\soclo$ holds for $m \geq k$.
The same holds in positive characteristic for the plus closure.

If $k \leq (d_1+d_2+d_3)/2$, then $\maxdeg_1(\Syz(k))=0$ and
$\mindeg_1(\Syz(k)) =(2k-d_1-d_2-d_3) \delta \leq 0$.

If moreover $k = (d_1+d_2+d_3)/2$, then the sheaf of syzygies is
strongly semistable.
\end{corollary}

\proof The primary syzygy yields an invertible quotient sheaf which
is isomorphic to $\det \Syz(k)$. For $k \geq \frac{d_1+d_2+d_3}{2}$
this quotient sheaf has degree $\geq 0$ and we are in the situation
of Proposition \ref{exactsequencecrit}(iv) and (v).

Now suppose $k \leq (d_1+d_2+d_3)/2$. We have $ \deg (\O_Y) = 0 \geq
k - (d_1+d_2+d_3)/2 = \mu (\Syz(k))$; hence we are in the situation
of Lemma \ref{exactsequenceranktwo}. If $k = (d_1+d_2+d_3)/2$, then
$\Syz(k)$ is the extension of the structure sheaf by itself; hence
its degree is $0$ and this follows again from Lemma
\ref{exactsequenceranktwo}. \qed

\begin{example}
\label{fermatexample} Consider a Fermat polynomial $x^k+y^k+z^k \in
K[x,y,z]$ and let $R=K[x,y,z]/(x^k+y^k+z^k)$. Let $f=x^{d_1}, \,
g=y^{d_2}, h=z^{d_3}$ such that $d_i \leq k$ and $d_1+d_2+d_3 \leq
2k$. Then $(x^{k-d_1},y^{k-d_2},z^{k-d_3})$ is a primary syzygy of
total degree $k$. Therefore by Corollary \ref{primarysemistable} we
get $R_{\geq k} \subseteq (x^{d_1}, y^{d_2}, z^{d_3})^\pasoclo$ and
we also get $R_{\geq k} \subseteq (x^{d_1}, y^{d_2}, z^{d_3})^{\gr}$
in positive characteristic.
\end{example}

\begin{corollary}
\label{primaryexclusion} Suppose the situation and notation of {\rm
\ref{situation}}. Suppose that there exists a primary syzygy of
total degree $k \leq (d_1+d_2+d_3)/2$. Suppose that the
characteristic of the algebraically closed field $K$ is zero or $p
\gg 0$. Let $f_0 \in R$ be a homogeneous element of degree $ \deg
(f_0)= m$. Then the following hold.

If $m <k$, then $f_0 \in (f_1,f_2,f_3)^\soclo $ if and only if $f_0 \in (f_1,f_2,f_3)$.

If $m \geq d_1+d_2+d_3-k$, then $ f_0 \in (f_1,f_2,f_3)^\soclo$.

\end{corollary}
\proof From Corollary \ref{primarysemistable} we get the short exact
sequence $0 \ra \O_Y(m-k) \ra \Syz(m) \ra \O_Y(m+k-d_1-d_2-d_3 ) \ra
0$, where $ \deg ( \O_Y(m-k)) \geq \mu (\Syz(m)) \geq
\deg(\O_Y(m+k-d_1-d_2-d_3))$. Thus we are in the situation of
Corollary \ref{exactsequencedecide}(i) and (iv). \qed

\medskip
We may also deduce a result about the plus closure.

\begin{corollary}
\label{exactsequence} Suppose the situation and notation of {\rm
\ref{situation}} and suppose that $K$ has positive characteristic.
Suppose that there exists a primary syzygy of total degree $k \leq
(d_1+d_2+d_3)/2$. Then $(f_1,f_2,f_3)^\soclo = (f_1,
f_2,f_3)^{\gr}$.
\end{corollary}

\proof Let $c \in H^1(Y,\Syz(m))$ denote the cohomology class of a
homogeneous element $f_0 \in R$ of degree $m$. We look at the
sequence from Corollary \ref{primarysemistable}, $ 0 \ra \O_Y(m-k)
\ra \Syz(m) \ra \O_Y(m+k-d_1-d_2-d_3) \ra 0$ and run through the
cases according to Corollary \ref{exactsequencedecide}. If $m \geq
d_1+d_2+d_3 -k$, then $f_0 \in (f_1, f_2,f_3)^+$ follows from
Proposition \ref{exactsequencecrit}(v). So suppose that $m <
d_1+d_2+d_3 -k $. If the image of $c$ in $H^1(Y,\O_Y(m+k
-d_1-d_2-d_3))$ is $\neq 0$, then $f_0 \not\in (f_1,f_2,f_3)^\soclo$
by \ref{exactsequencecrit}(i). So we may assume that $c$ stems from
$e \in H^1(Y, \O_Y(m-k))$. If $m \geq k$ or $c=0$, then the
pull-back of $e$ under a finite mapping is zero; hence $f_0 \in
(f_1, f_2,f_3)^+$. If $m < k$ and $c \neq 0$, then $f_0 \not\in
(f_1, f_2,f_3)^\soclo$ by Proposition \ref{exactsequencecrit}(ii).
\qed

\begin{example}
Consider the ideal $(x^{10}, y^{10}, z^{10})$ on the curve given by
the equation $x^4+y^4 =z^4$. Due to Remark \ref{riemannroch} there
exists a syzygy $\neq 0$ of degree $16$. We have
$$z^{16} = (x^4+y^4)^4= x^{16}+ 4 x^{12}y^4+6x^8y^8 + 4x^4y^{12}+y^{16}$$
and
$$ x^4z^{12} =x^4 (x^4+y^4)^3= x^{16}+ 3 x^{12}y^4+3x^8y^8 + x^4y^{12} \, . $$
Therefore we may write
$$z^{10} (z^6-2z^2x^4) =
-x^{16} -2x^{12}y^4 +2x^4y^{12} +y^{16} =x^{10}(-x^6-2x^2y^4)
+y^{10} (2x^4y^2+y^6) \, ;$$ hence we have the syzygy $(z^6-2z^2x^4,
x^6+2x^2y^4, -2x^4y^2-y^6)$ of total degree $16$. This syzygy is
primary: if $x=0$ or $y=0$, then  $x=y=z=0$, so suppose $x ,y \neq
0$. Then $x^4 +2y^4 =0$ and $2x^4+y^4=0$ which gives $3x^4=0$; hence
the syzygy has no common zero (in characteristic $\neq 3$) and is
therefore primary. We have therefore the short exact sequence
$$0 \lra \O_Y \lra \Syz(16) \lra \O_Y(2) \lra 0 \, $$
and it follows from Corollary \ref{primarysemistable}
that $R_{\geq 16} \subset (x^{10},y^{10},z^{10})^\soclo$.
We tensor this short exact sequence with $\O(k-16)$ and dualize it to get
$$0 \lra \O_Y(14-k) \lra \Syz^\dual( -k) \lra \O_Y(16-k) \lra 0 \, .$$
This shows that $\Syz^\dual (-k)$ is ample for $k \leq 13$ as an
extension of two ample invertible sheaves, therefore $(x^{10},
y^{10}, z^{10})^\soclo \cap R_{\leq 13} = (x^{10}, y^{10}, z^{10})
\cap R_{\leq 13}$ holds in characteristic $0$ and $p \gg 0$. If
there does not exist a syzygy of degree $15$, then this also holds
for $k=14$ due to Proposition \ref{amplecritrelation}.
\end{example}

\begin{example}
\label{examplerelationindecomposable }
Let $K$ denote an algebraically closed field
and consider
$$R=K[x,y,z]/(x^\delta+ay^\delta+bz^\delta+cxz^{\delta-1}+dyz^{\delta-1}) \,  $$
where $a,b,c,d \neq 0 \, $ are chosen such that $Y= \Proj \, R$ is
smooth. Consider the syzygy sheaf for the elements $x^\delta,\,
y^\delta,\, z^\delta$. Then we have a syzygy of total degree $\delta
+1$, given by $(z,az, bz+cx+dy)$. Since $\delta +1 < \frac{3 \delta
}{2} - \frac{\delta-3}{2} = \frac{3 \delta}{2} - \frac{g-1}{\delta}
$, it follows from Proposition \ref{maxinvertbound}(iii) that $\Syz$
is decomposable.

The syzygy $(z,az, bz+cx+dy)$ is primary if and only if $cx+dy$ and
$x^\delta + ay^\delta$ have no common homogeneous zero. This is true
if and only if $(-\frac{d}{c})^\delta \neq -a$. If this is true,
then we have the splitting $\Syz( \delta +1 )= \O_Y \oplus \O_Y(-
\delta +2 )$, where the second summand corresponds to a syzygy of
total degree $2 \delta -1$. We can find such a syzygy in the
following way: There exists a polynomial $P(x,y)$ in $x$ and $y$ of
degree $\delta -1$ such that $(cx+dy) P(x,y) = rx^\delta
+sy^\delta$. Then $(P+rz^{\delta -1},aP+sz^{\delta -1}, bP )$ is a
syzygy of total degree $2\delta-1$, since $Px^\delta + aPy^\delta +
rx^\delta z^{\delta -1} + sy^\delta z^{\delta -1} +  bPz^\delta =
Px^\delta +aPy^\delta +Pbz^\delta +Pcxz^{\delta -1} +Pdyz^{\delta
-1}=0$.
\end{example}

\begin{corollary}
\label{bastel} Let $f_1,f_2,f_3 \in K[x,y,z]$ be homogeneous
polynomials of degree $d_1,d_2,d_3$ such that $d_1+d_2+d_3=2k$ is
even and $k \geq d_i$ for $i=1,2,3$. Let $g_1,g_2,g_3 \in K[x,y,z]$
be homogeneous of degree $k-d_i$. Suppose that
$V(f_1,f_2,f_3)=V(g_1,g_2,g_3)=V(x,y,z)$. Set
$F=f_1g_1+f_2g_2+f_3g_3$ and suppose that $R=K[x,y,z]/(F)$ is a
normal domain. Then the sheaf of syzygies $\Syz(m)$ for
$f_1,f_2,f_3$ on $Y= \Proj R$ is an extension of the structure sheaf
by itself and is strongly semistable. In particular
$$(f_1,f_2,f_3)^\soclo = R_{\geq k} +(f_1,f_2,f_3) \, .$$
If furthermore the characteristic of $K$ is positive,
then $(f_1,f_2,f_3)^\soclo = (f_1,f_2,f_3)^{+ \rm gr}$.
\end{corollary}

\proof The syzygy $(g_1,g_2,g_3)$ is primary of total degree $k$;
thus this follows from Corollary \ref{primarysemistable}. \qed

\section{The tight closure of $(x^{a},y^{a},z^{a})$ in $K[x,y,z]/(F)$}
\label{examplea}

In this section we study the tight closure of $(x^{a},y^{a},z^{a})$
in $R=K[x,y,z]/(F)$, where $F$ is an  irreducible homogeneous
polynomial of degree $\delta $ such that $R$ is normal. The expected
generic answer is by Theorem \ref{slopetightcriteria}(iii) that $
(x^a,y^a,z^a)^\soclo= (x^a,y^a,z^a) + R_{\geq \frac{3}{2}a}$. We
have however to check that the sheaf of syzygies is semistable to
obtain this result in characteristic zero, and in positive
characteristic we have to do even more. Sullivant has made some
computations implemented in Macaulay2 for the monomial ideals
$(x^{a},y^{a},z^{a})^\soclo$ for $a =2,3,4$ for small prime numbers
$p$ and small degree $\delta$ for the Fermat equations $x^\delta
+y^\delta +z^\delta=0$ (see \cite{sullivant}). These computations
have led him to conjectures about the behavior of the tight closure;
we will prove some of his conjectures.

The homogeneous ideal generators $x^a,y^a,z^a$ yield the sheaf of
syzygies $\Syz(m) = Rel (x^a,y^a,z^a)(m)$ on the smooth projective
curve $Y= \Proj R \subset \PP^2$. This sheaf is the restriction of
the sheaf of syzygies $\Syz_{\PP^2}(m)=Rel_{\PP^2} (x^a,y^a,z^a)(m)$
on the projective plane. On $\PP^2$ we have the presenting sequence
$$0 \lra \Syz_{\PP^2}(m) \lra \bigoplus_3 \O_{\PP^2} (m - a)
\stackrel{x^a,y^a,z^a}{\lra} \O_{\PP^2}(m) \lra 0 \, $$
and the exact sequence (from the Koszul complex)
$$0 \lra \O_{\PP^2}(m- 3a)  \lra
\bigoplus_3 \O_{\PP^2}(m-2a) \lra \Syz_{\PP^2}(m) \lra 0 \, ,$$
where the surjection is given by the standard syzygies
$$(-y^a,x^a,0), \hspace{15pt} (z^a,0,-x^a),\hspace{15pt} (0,-z^a,y^a) \, $$
and the injection by $1 \mapsto (z^a,y^a,x^a)$.
Since the sheaves in these short exact sequences are locally free,
their restrictions to a curve are also exact (they are subbundles).

\begin{lemma}
\label{noglobalsections} Let $F \in K[x,y,z]$ denote an irreducible
homogeneous polynomial of degree $\delta $ and suppose that it
defines a smooth projective curve $Y$. Let $\Syz(m)$ denote the
sheaf of syzygies on $Y$ for the elements $x^a,y^a,z^a$, $a \geq 1$.
Then every syzygy $ \in \Gamma(Y, \Syz(k))$ for $ k < \delta$ is a
linear combination of the three standard syzygies $(-y^a,x^a,0)$,
$(z^a,0,-x^a)$ and $(0,-z^a,y^a)$. In particular $ \Gamma(Y,
\Syz(k))=0 $ for $k < 2a, \delta$.
\end{lemma}
\proof On $\PP^2$ we have the exact sequence
$$0 \lra \O_{\PP^2}(-a)  \lra \bigoplus_3 \O_{\PP^2}
\lra \Syz_{\PP^2}(2a) \lra 0 \, ,$$ where the surjection is given by
the standard syzygies. This surjection is also globally a
surjection. We tensor the exact sequence for $Y \subset \PP^2$, that
is,
$$ 0 \lra \O_{\PP^2}(- \delta) \stackrel{F}{\lra}
\O_{\PP^2} \lra \O_Y \lra 0 \, ,$$ with $\Syz_{\PP^2}(k)$, and by
applying $\Gamma( \PP^2, - )$ we get
$$ \Gamma(\PP^2,\Syz_{\PP^2}(k)) \lra  \Gamma(\PP^2, \Syz_Y (k)) \lra
H^1(\PP^2,\Syz_{\PP^2}(k- \delta )) \, .$$ We want to show that the
term on the right is zero. By the presenting sequence for the
syzygies we get
$$0 \lra \Syz_{\PP^2}(k-\delta) \lra \bigoplus_3 \O_{\PP^2} (k-\delta -a)
\lra \O_{\PP^2}(k - \delta) \lra 0 \, $$
and hence
$$ \Gamma(\PP^2, \O_{\PP^2}(k - \delta)) \lra
H^1(\PP^2,\Syz_{\PP^2}(k- \delta )) \lra \bigoplus_3 H^1(\PP^2,
\O_{\PP^2}(k- \delta-a))=0 \, $$ shows that
$\Gamma(\PP^2,\Syz_{\PP^2}(k)) \lra  \Gamma(\PP^2, \Syz_Y (k))$ is
surjective for $k < \delta$. Then the statements follow since they
are true on $\PP^2$. \qed

\begin{proposition}
\label{semistablea} Let $F \in K[x,y,z]$ denote a polynomial of
degree $\delta \geq 3a-1 $ and suppose that it defines a smooth
projective curve $Y$. Then the sheaf of syzygies $\Syz(m)$ for the
elements $x^a,y^a,z^a$ is semistable.
\end{proposition}

\proof Suppose that there exists an invertible subsheaf $\shL
\subset \Syz( \delta -1)$ of degree $ > \mu (\Syz( \delta -1)) =
\frac{2( \delta -1)-3 a}{2} \delta =  (\delta -1-\frac{3}{2}a)
\delta $. Then
\begin{eqnarray*}
\deg (\shL) + 1 - g   &> &
(\delta - 1 -\frac{3}{2}a) \delta + 1 - \frac{(\delta -1) (\delta -2)}{2} \cr
&=& \delta ( \delta - 1 -\frac{3}{2}a - \frac{\delta -3}{2} )
= \delta ( \frac{\delta}{2} - \frac{3}{2}a + \frac{1}{2})
\, .
\end{eqnarray*}
This expression is $\geq 0$ for $ \delta \geq 3a-1$; hence due to
the theorem of Riemann-Roch the invertible sheaf $\shL$ has global
sections $\neq 0$. Therefore we look at the global sections of
$\Syz( \delta -1)$ and study their zeros in order to get a bound for
the maximal degree of a subbundle.

We use the inclusion $ \Syz(\delta -1 ) \subset \O_Y( \delta -a-1)
\oplus \O_Y(\delta -a-1) \oplus \O_Y(\delta -a -1) $ and think of a
global syzygy $S \in \Gamma(Y, \Syz(\delta -1))$ as given by three
polynomials $(S_1, S_2,S_3)$ of degree $ \delta -a-1$. By Lemma
\ref{noglobalsections} we know that the global syzygies of total
degree $\delta -1$ are of the form
$$ (S_1,S_2,S_3) =  A(-y^{a},x^{a},0) +B(-z^{a},0,x^{a}) +C(0,z^{a},-y^{a}) \, ,$$
where $A,B,C$ are homogeneous polynomials in $K[x,y,z]$ of degree $
\delta -2a -1$. A syzygy $ \O_Y \ra \Syz(\delta -1)$ has a zero if
and only if the three components $S_1=Ay^{a} + Bz^{a}$,
$S_2=Ax^{a}+Cz^{a}$ and $S_3=Bx^{a}-Cy^{a}$ have a common zero. We
have to show that the number of zeros (with multiplicities) of such
a syzygy is bounded by the slope of $\Syz( \delta -1)$, that is, by
$ ( \delta -1 - \frac{3}{2}a )  \delta $. We write
$$ S_1= Ay^{a}+ Bz^{a} =Q_1 P \hspace{15pt} \mbox{   and } \hspace{15pt}  S_2= Ax^{a}+Cz^{a} = Q_2 P \, ,$$
where $Q_1$, $Q_2$ have no common divisor, $\deg (Q_1) = \deg (Q_2)
= t$. We assume first that $ a \leq t \leq \delta -a-1$. The number
of zeros of the syzygy on the curve given by $F=0$ is then bounded
by ($F$ and $P$ have no common divisor, since $F$ is irreducible of
degree $\delta > \deg (P)$)
$$ \deg (P) \delta + \deg (Q_1) \deg (Q_2) =
(\delta - a-1 -t) \delta + t^2 \, .$$
We claim that
$ t^2 \leq \delta ( t - \frac{1}{2}a) $.
We can check this at the boundaries
$t=a$ and $t= \delta -a-1$.
For $t=a$ this is the inequality $a^2 \leq \delta a/2$,
which is true since $\delta \geq 3a-1 \geq 2a$ for $a \geq 1$.
For $t= \delta -a-1$
we have to show that
$ (\delta -a-1)( \delta -a-1) \leq \delta (\delta - \frac{3}{2}a -1)$,
which yields the condition
$a^2+2a+1 \leq \delta (1 + a/2)$.
But this is true since $\delta \geq 3a-1$ and $\delta \geq 1$.

The claim implies that
$ -\delta (a+t) + t^2 \leq - \frac{3}{2}a \delta$.
Therefore we have
$$(\delta - a-1 - t ) \delta + t^2 = (\delta -1) \delta
- \delta (a +t) + t^2
\leq (\delta -1) \delta  -\frac{3}{2}a \delta \, .$$

Assume now that $t < a$.
We may assume that the powers $z^b$, $b \geq a$, do not occur in the polynomial $A$,
since $z^{a}(-y^{a},x^{a},0) = y^{a}(-z^{a},0,x^{a}) + x^{a}(0,z^{a},-y^{a})$.
{}From $Q_2S_1=Q_1S_2$ we obtain the equation
$$A(Q_1x^{a}-Q_2y^{a}) = z^{a}(Q_2B-Q_1C) \, .$$
If $z$ divides $Q_1x^{a}-Q_2y^{a}$, then it would also (since $\deg
(Q_1)=t < a$) divide $Q_1$ and $Q_2$, but they are coprime. Hence
$z^a$ divides $A$, but then $A=0$. Therefore $Q_2B=Q_1C$ and hence
$B= Q_1 D$ and $C=Q_2D$, where $D$ is a polynomial of degree $
\delta -2a-1-t $. A zero of the syzygy is given by $DQ_1z^a=0$,
$DQ_2z^a=0$ and $DQ_1x^a - D Q_2y^a=0$, hence by $D=0$ or by
$Q_1z^a=Q_2z^a=Q_1x^a - Q_2y^a=0$. The polynomials $Q_1z^a$ and
$Q_1x^a - Q_2y^a$ do not have a common divisor, since $Q_1$ and
$Q_2$ are coprime and since $z$ does not divide $Q_1x^{a}-Q_2y^a$.
Therefore the number of zeros is bounded by $( \delta-2a-1-t) \delta
+ (t+a)^2$. Again we have to check that this is $\leq   ( \delta -1
- \frac{3}{2}a )  \delta $ and this is equivalent to $(t+a)^2 \leq
\delta (\frac{1}{2}a +t)$. But this is true for $t=a-1$ and $t=0$.
\qed

\begin{corollary}
\label{socloa}
Let $K$ denote an algebraically closed field of characteristic zero.
Let $F \in K[x,y,z]$ denote an irreducible homogeneous polynomial
of degree $\delta $ and suppose that $R=K[x,y,z]/(F)$ is normal.
Then for $\delta \geq 3a-1$ we have
$$ (x^{a},y^{a},z^{a})^\soclo =
(x^{a},y^{a},z^{a}) + R_{\geq \frac{3}{2}a} \, .$$
\end{corollary}
\proof This follows from Proposition \ref{semistablea} and Theorem
\ref{slopetightcriteria}(iii). \qed

\medskip
We also obtain from \ref{semistablea} results in positive
characteristic. We deduce first the following ampleness result.

\begin{corollary}
\label{amplea} Let $F \in K[x,y,z]$ denote an irreducible
homogeneous polynomial of degree $\delta$ such that $R=K[x,y,z]/(F)$
is normal. Let $\Syz(m)$ denote the sheaf of syzygies for the
elements $x^{a},y^{a},z^{a}$ of total degree $m$ on $Y= \Proj \, R$.
Suppose that $ \delta \geq 3a-1$ and that the characteristic of $K$
is zero or $p \geq \delta -3$. Then $\Syz(m)$ is an ample sheaf for
$m > \frac{3}{2}a$.
\end{corollary}

\proof We will use the ampleness criterion of Hartshorne-Mumford for
locally free sheaves of rank two on a smooth projective curve $Y$
(see \cite[Proposition 7.5 and Corollary 7.7]{hartshorneample}). It
states that $\shS$ is ample provided that

\renewcommand{\labelenumi}{(\roman{enumi})}
\begin{enumerate}

\item
$ \deg (\shS) > \frac{2}{p} (g-1)$, where $g$ is the genus of $Y$,

\item
every invertible quotient sheaf $ \shS \ra \shL \ra 0 $ has positive
degree,
\end{enumerate}
(in characteristic zero the first condition is just that $\deg
(\shS) >0$). We have $\deg (\Syz(m)) = (2m- 3a) \delta > \delta$ and
$g-1= \delta (\delta -3) /2$; therefore the first condition is
fulfilled for $p \geq \delta -3$. The second condition follows from
Proposition \ref{semistablea}: by semistability we have $\deg (\shL)
\geq \mu (\Syz(m)) >0$. \qed

\begin{corollary}
\label{frobeniustight}
Let $K$ denote an algebraically closed
field of positive characteristic $p$ and let
$F \in K[x,y,z]$ denote a homogeneous polynomial of degree $\delta$ such that
$R=K[x,y,z]/(F)$ is a normal domain.
Suppose that $ \delta \geq 3a-1$.
Then the following hold for the Frobenius closure and the tight closure of
$(x^{a},y^{a},z^{a})$.

\renewcommand{\labelenumi}{(\roman{enumi})}
\begin{enumerate}

\item
$R_m \subset (x^{a},y^{a},z^{a})^F $ for $m > \frac{3}{2} a$ and
$p \geq \delta -3 $.

\item
$R_m \subset (x^{a},y^{a},z^{a})^\soclo $ for $m > \frac{3}{2} a$ and
$p \geq \delta -3 $.

\item
$(x^{a},y^{a},z^{a})^\soclo \cap R_m = (x^{a},y^{a},z^{a}) \cap R_{m}$
for $m < \frac{3}{2} a$ and $ p \gg 0$.

\end{enumerate}

\end{corollary}

\proof (i) We know by Corollary \ref{amplea} that the sheaf of
syzygies $\Syz(m)$ is ample for $m > \frac{3}{2}a$ and $p \geq
\delta -3 $. Let $c \in H^1(Y, \Syz(m))$ denote a cohomology class
(given by an element $f_0 \in R_m$). An ample sheaf $\shS$ on a
smooth projective curve over an algebraically closed field of
positive characteristic is also cohomologically $p$-ample. This
means that for every coherent sheaf $\shF$ we have $H^{i}(Y,
\shS^{(q)} \otimes \shF)=0$ for $i \geq 1$ and $q \gg 0$, where
$\shS^{(q)}$ denotes the Frobenius pull-back of $\shS$ (see
\cite{hartshorneample}, \cite{migliorinicohomological}). In
particular the mapping $H^1(Y, \Syz(m)) \ra H^1(Y, q^*\Syz(m))
=H^1(Y, \Syz(x^{aq},y^{aq},z^{aq})(qm) )$ is $0$ for $q \gg 0$.
Hence $c^q =0$ and therefore $f^q \in (x^{aq},y^{aq},z^{aq})$. This
means that $f$ belongs to the Frobenius closure of $(x^a, y^a,z^a)
$.

(ii) This follows from (i).

(iii) We have $\Syz(m)^\dual = \Syz(m) \otimes \O_Y(3a-2m) = \Syz(
3a-m)$, so for $m < \frac{3}{2}a$ the dual sheaf $\shF(-m) =
\Syz(m)^\dual$ is ample. A cohomology class $c \in H^1(Y, \Syz(m))$
defines an extension $0 \ra \O_Y \ra \shF' \ra \shF(-m) \ra 0$. In
characteristic zero and for $c \neq 0$ this extension $\shF'$ is
also an ample sheaf. This is then also true for $p \gg 0$; hence the
complement $\PP(\shF') - \PP(\shF(-m))$ is an affine scheme for $p
\gg 0$. This means that $f \not \in (x^{a},y^{a}, z^{a})^\soclo$ for
$ p \gg 0$ with the exception that $c= \delta (f) =0$, which means
that $f \in (x^{a},y^{a}, z^{a})$. \qed

\begin{remark}
Note the difference between $a$ even and $a$ odd. For $a$ odd the
critical value $\frac{3}{2}a$ is not the degree of a polynomial and
so Corollary \ref{frobeniustight} gives a complete answer for the
tight closure of $(x^{a},y^{a}, z^{a})$. This also explains that
Sullivant's algorithm works much better for degree $3$ than for
degree $2$ or $4$.
\end{remark}

\begin{remark}
The result in Corollary \ref{frobeniustight}(i) is under the much
stronger condition $m > 2a$ easy to prove using the Koszul
resolution. This yields the surjection
$$\bigoplus_3 H^1(Y, \O_Y((m-2a)q))
\ra H^1(Y, \Syz(x^{aq},y^{aq},z^{aq})(mq)) \ra 0$$ and $h^1
(\O_Y((m-2a)q)) = h^0 (\O_Y( (2a-m)q) \otimes \O_Y(\delta -3)) =0$
for $q \gg 0$.
\end{remark}

\begin{remark}
The results above imply in particular Sullivant's conjectures 4.2
and 4.5 about the behavior of $(x^3,y^3,z^3)^\soclo$ ($a=3$) in the
Fermat rings for $ \delta \geq 8$ and $ p \gg 0$. For degree $
\delta \geq 4$ see Section \ref{tightexample}.
\end{remark}

We close this section with an easy observation concerning the strong
semistability of the sheaf of syzygies.

\begin{proposition}
\label{nonstrongly} Let $\Syz(m)$ denote the sheaf of syzygies for
the elements $x^a,y^a$, $z^a$ on the Fermat curve given by $x^\delta
+ y^\delta +z^\delta =0$ over an algebraically closed field of
positive characteristic $p$. Suppose that $\frac{3a}{2}q > \delta q'
\geq aq$ for prime powers $q$ and $q'$ of $p$. Then $\Syz$ is not
strongly semistable.
\end{proposition}

\proof The pull-back of $\Syz(x^a,y^a,z^a)(0)$ under the $k$-th
Frobenius morphism, $q= p^k$, is $\Syz(x^{aq},y^{aq},z^{aq})(0)$.
Since $\delta q' \geq aq $, we can use the $q'$-power of the curve
equation, that is, $x^{\delta q'} +y^{\delta q'}+ z^{\delta q'} =0$
directly as the syzygy $ (x^{\delta q'-aq},\! y^{\delta q'-aq},\!
z^{\delta q'- aq})$ for $x^{aq}, y^{aq}, z^{aq}$ of total degree
$\delta q'$. This gives a non-trivial morphism $\O_Y \ra
\Syz(x^{aq},y^{aq},z^{aq})( \delta q')$. The slope on the right is $
\delta q'-\frac{3}{2} {aq} <0 $, so this sheaf of syzygies is not
semistable. \qed

\begin{remark}
The condition in Proposition \ref{nonstrongly}
means either $\frac{3a}{2}p^k > \delta \geq ap^k$ or
$\frac{3a}{2} > \delta p^k \geq a$.
For $\delta p^k \geq \frac{3}{2} a$ we get from a similar argument
that $ R_{\geq \delta p^k} \subset (x^a,y^a,z^a)^\soclo$
(and even in the plus closure).
\end{remark}

\section{The tight closure of $(x^2,y^2,z^2)$ and $(x^3,y^3,z^3)$}
\label{tightexample}

We continue with the computation of $(x^a,y^a,z^a)^\soclo$ in
$R=K[x,y,z]/(F)$, where $F$ is an irreducible homogeneous polynomial
of degree $\delta$ such that $R$ is normal. We deal first with the
case $a=2$ and consider also the behavior for small degrees and
small prime numbers $p$. The computation of $(x^2,y^2,z^2)^\soclo$
is in fact a question of whether $xyz \in (x^2,y^2,z^2)^\soclo$
holds or not.

\medskip
If the degree of $F$ is one, then we may write $x=ay+bz$ and hence
$xyz =ay^2z + byz^2 \in (x^2,y^2,z^2)$, so $xyz$ belongs to the
ideal itself. Suppose that the degree of $F$ is two; hence $F$
defines a quadric. Suppose first that at least one of the mixed
monomials $xy$, $xz$, $yz$ occurs in $F$, say $xy$. Then the
multiple $zF$ shows again that $xyz \in (x^2,y^2,z^2)$. If however
$F=ax^2+by^2+cz^2$, then $K[x,y,z]/(F,x^2,y^2,z^2) =
K[x,y,z]/(x^2,y^2,z^2)$ and $xyz \not\in (x^2,y^2,z^2)$. If $F$ is
irreducible ($a,b,c \neq 0$), then $\Spec \, R$ is a normal cone
over the projective line; hence it is $F$-regular and $xyz \not\in
(x^2,y^2,z^2)^\soclo$. If $F$ is the product of two linear forms,
then $xyz \in (x^2,y^2,z^2)^\soclo$, since this (even $\in
(x^2,y^2,z^2)$) is true on the two planes and the containment to the
solid closure may be checked on the components.

\medskip
Now suppose that $F \in K[x,y,z]$ has degree $3$ and defines an
elliptic curve $Y$. If the coefficient of $F$ in $xyz$ is not zero,
then of course $xyz \in (x^2,y^2,z^2)$. Thus we may write
$F=Sx^2+Ty^2+Uz^2$, so that $(S,T,U)$ is a homogeneous syzygy for
$(x^2,y^2,z^2)$ of total degree $3$. If this syzygy is primary,
i.e., $V(S,T,U)=V(R_+)$, then we have a short exact sequence
$$0 \lra \O_Y \lra \Syz(3) \lra \O_Y \lra 0$$
and therefore $\Syz$ is semistable on the elliptic curve $Y$. Hence
$xyz \in (x^2,y^2,z^2)^\soclo$.

If however the syzygy $(S,T,U)$ is not primary, e.g., for
$F=x^3+y^3+(x+ay)z^2$, then we have a decomposition $\Syz(3) = \O(P)
\oplus \O(-P)$ ($P$ a point) and $xyz \not\in (x^2,y^2,z^2)^\soclo$,
since $H^1(Y,\Syz(3)) = H^1(Y, \O(-P))$ and by Proposition
\ref{exactsequencecrit}(2).

\medskip
We now examine the case where the equation of the curve is given
by a polynomial of degree $4$.
Note that Proposition \ref{semistablea} gives the semistability
only for $\delta \geq 5$.

\begin{lemma}
\label{degreefoursemistable} Let $F \in K[x,y,z]$ denote an
irreducible homogeneous polynomial of degree $4$ and suppose that it
defines a smooth projective curve $Y$. Let $\Syz(m)$ denote the
sheaf of syzygies on $Y$ for the elements $x^2,y^2,z^2$. Then
$\Syz(m)$ is semistable.
\end{lemma}

\proof Suppose that $\shL \ra \Syz(3)$ is non-trivial with $\deg
(\shL) >0$. Then $\shL \otimes \O_Y(1) \ra \Syz(4)$ is non-trivial
and $\deg (\shL(1)) \geq 5$; hence $\shL(1)$ has global sections
$\neq 0$ due to the theorem of Riemann-Roch. It is then enough to
show that every non-trivial global section $\O_Y \ra \Syz(4)$ has at
most four zeros (counted with multiplicities). Such a section is
given by a syzygy $S=(S_1,S_2,S_3)$, where the $S_i$ are homogeneous
polynomials of degree $2$. Since $T=S_1x^2+S_2y^2+S_3z^2=0$ on the
curve, we have $T= \lambda F$ in $K[x,y,z]$, $\lambda \in K$.

The zeros of the section $S$ are the common zeros of $(S_1,S_2,S_3)$.
If two of the $S_1,S_2,S_3$ have no common divisor, say $S_1$ and $S_2$,
then $V_+(S_1) \cap V_+(S_2)$ consists of four points
(counted with multiplicities).

So suppose that $S_1,S_2,S_3$ have together one common linear factor
$P$. Then we may write $P(Q_1x^2+Q_2y^2+Q_3z^2)= \lambda F$ in
$K[x,y,z]$. For $\lambda \neq 0$ the polynomial $F$ would be
reducible, which is excluded. Hence $\lambda =0$. Here $P=0$ would
imply that we are dealing with the zero syzygy; hence
$(Q_1,Q_2,Q_3)$ is a global syzygy of degree $3$. But these are all
trivial due to Lemma \ref{noglobalsections}; hence this case is not
possible.

So suppose that the $S_i$ have pairwise one common linear factor.
Then the common zeros of $S_1,S_2,S_3$ are the three
intersection points of the triangle.
\qed

\begin{remark}
The global sections of the sheaf of syzygies $\Syz(x^2,y^2,z^2)(4)$
on a curve of degree $4$ depend heavily on the curve equation $F=0$.
We always have the standard syzygies and their linear combinations,
but every way of writing $F \in (x^2,y^2,z^2)$ yields also a syzygy.
\end{remark}

\begin{corollary}
Let $K$ denote an algebraically closed field of characteristic zero.
Let $F \in K[x,y,z]$
denote an irreducible homogeneous
polynomial of degree $\delta =4 $ and suppose that
$R=K[x,y,z]/(F)$ is normal.
Then $xyz \in (x^2,y^2,z^2)^\soclo$.
\end{corollary}
\proof Lemma \ref{degreefoursemistable} shows that the sheaf of
syzygies is semistable; hence the numerical criterion Theorem
\ref{slopetightcriteria}(iii) shows that $R_{\geq 3} \subseteq
(x^2,y^2,z^2)^\soclo$. \qed

\begin{example}
\label{fermat4} We have a closer look at the Fermat quartic
$F=x^4+y^4+z^4$, $\Char (K) \neq 2$. We consider the syzygy for
$x^2,y^2,z^2$ of total degree $4$ given by
$$(-y^2-z^2, x^2+ \sqrt{2}i z^2, x^2- \sqrt{2} i y^2)
=(-y^2,x^2,0)+ (-z^2,0,x^2) +\sqrt{2} i(0,z^2, -y^2) \, .$$ For
$z=0$ this syzygy has no zero on the curve. For $z \neq 0$ we find
the four zeros
$$z=1,\hspace{15pt} y= +/- i, \hspace{15pt} x = +/- \sqrt[4]{2} \sqrt{-i} $$
and no more. Denote them by $P_1, P_2,P_3, P_4$ and let $\Sigma =
P_1+ P_2+P_3+ P_4$ be their Weil divisor. The syzygy $\O_Y \ra
\Syz(4)$ then factors through $\O_Y \ra \O_Y( \Sigma) \ra \Syz(4)$
and we get an invertible subsheaf of $\Syz(4)$ without zero, hence a
subbundle. This yields a short exact sequence
$$ 0 \lra \O_Y( \Sigma ) \otimes \O_Y(-1) \lra \Syz(3)
\lra \O_Y(-\Sigma ) \otimes \O_Y(1) \lra 0 \, .$$ The degree on the
left and on the right is $0$. {}From this it follows not only that
$\Syz(m)$ is semistable, but also that it is strongly semistable by
Lemma \ref{exactsequenceranktwo}. Hence $xyz \in
(x^2,y^2,z^2)^\soclo$ holds on $x^4+y^4+z^4=0$ in positive
characteristic $p \geq 3$ also.
\end{example}

\begin{remark}
Note the difference between the sequence in degree three and in the
Fermat example of degree four. Both show that $\Syz$ is strongly
semistable but not stable. The first sequence shows at once that
$xyz \in (x^2,y^2,z^2)^+$, which is not clear at all in Example
\ref{fermat4}.
\end{remark}

\begin{example}
Let $F=zx^3+xy^3+yz^3=0$. If we consider this equation as a syzygy
for $x^2,y^2,z^2$, then this has exactly three zeros (as in the last
part of the proof of Lemma \ref{degreefoursemistable}). We get a
sequence
$$0 \lra \O(P_1+P_2+P_3 -H) \lra  \Syz(3) \lra \O(H-P_1-P_2-P_3)\lra 0 \, .$$
Is $\Syz(3)$ stable? Is it strongly semistable in positive
characteristic?
\end{example}

We look now at the situation of $\deg (F)= \delta \geq 5$. We know
that $xyz \in (x^2,y^2,z^2)^\soclo$ holds in $R=K[x,y,z]/(F)$ for
$\deg (F) = \delta \geq 5$ in characteristic zero due to Corollary
\ref{socloa}. This is in general not true in positive
characteristic. It is not even clear whether or not $xyz \in
(x^2,y^2,z^2)^\soclo $ holds on $x^\delta + y^\delta +z^\delta$,
$\delta \geq 5$, for infinitely many or almost all prime
characteristics.

\begin{example}
We consider the ideal $(x^2, y^2,z^2)$ on the curve given by
$x^7+y^7+z^7=0$ for characteristic $p=3$. The curve equation gives
at once a global syzygy for the elements $x^6,y^6,z^6$ (the third
power of $x^2, y^2,z^2$) of total degree $7$. Therefore we have the
short exact sequence $0 \ra \O_Y \ra \Syz(x^6,y^6,z^6)(7) \ra
\O_Y(-4) \ra 0$ showing that $\Syz(x^2,y^2,z^2)$ is not strongly
semistable (see also Proposition \ref{nonstrongly}). To decide
whether $xyz \in (x^2, y^2,z^2)^\soclo$ holds, we have to look at
$$0 \lra \O_Y(2) \lra \Syz(x^6,y^6,z^6)(9) \lra \O_Y(-2) \lra 0 \, .$$
The element $(xyz)^3$ yields the cohomology class
$$( \frac{(xyz)^3}{x^6},- \frac{(xyz)^3}{y^6} ,0)
\in H^1(Y, \Syz(x^6,y^6,z^6)(9)) \, ,$$ which maps to $ \frac{
(xyz)^3 z}{x^6 y^6} = \frac{z^4}{x^3y^3} \in H^1(Y, \O_Y(-2))$ by
Remark \ref{dual}. This class is not zero. Therefore the (dual)
extension $0 \ra \O_Y \ra \shG \ra \O_Y(2) \ra 0$ given by this
class is not trivial. {}From this and from $ \deg (\shG)=14 >
\frac{2}{3} (15-1)$ we see that $\shG$ is ample. Hence $z^4 \not\in
(x^3,y^3)^\soclo$ and then $xyz \not\in (x^2, y^2,z^2)^\soclo$ also.
\end{example}

We now look at $(x^3,y^3,z^3)^\soclo$ in order to extend the results
from the last section for the Fermat equations $x^\delta +y^\delta
+z^\delta =0$ of low degrees. For $\delta =1,2$ we have
$(x^3,y^3,z^3)^\soclo =(x^3,y^3,z^3)$. For $\delta =3$ we have
$(x^3,y^3,z^3)^\soclo = (x^3,y^3)^\soclo = (x^3,y^3) + R_{\geq 6}$.

\begin{example}
For $x^4+y^4+z^4=0$ we have the short exact sequence $0 \ra \O_Y \ra
\Syz(4) \ra \O_Y(-1) \ra 0$ given by the syzygy $(x,y,z)$. Therefore
the sheaf of syzygies $\Syz(m)$ for $x^3,y^3,z^3$ is not semistable
on the Fermat quartic. The sheaf $\Syz(5)$ is also not ample. It is
however the extension of two invertible sheaves of degree $\geq 0$;
hence it follows by Proposition \ref{exactsequencecrit} that
$R_{\geq 5} \subset (x^3,y^3,z^3)^\soclo$ and also that $R_{\geq 5}
\subset (x^3,y^3,z^3)^+$ in positive characteristic.
\end{example}

\begin{example}
\label{fermatfuenf} For $x^5+y^5+z^5=0$ we have the short exact
sequence $0 \ra \O_Y \ra \Syz(5) \ra \O_Y(1) \ra 0$ given by the
syzygy $(x^2,y^2,z^2)$. This extension is not trivial, since there
does not exist a non-trivial mapping $\O_Y \ra \Syz(4)$. Therefore
the sheaf $\Syz(5)$ is ample in characteristic zero and hence also
for $p \gg 0$. To obtain a bound for the prime number, we use the
criterion of Mumford-Hartshorne. Suppose that there exists an
invertible quotient sheaf $\Syz(5) \ra \shL \ra 0$ of $\deg (\shL)
\leq 0$. The exact sequence shows at once that $\shL = \O_Y$ and
hence that it would split. Hence it follows that
$(x^3,y^3,z^3)^\soclo =(x^3,y^3,z^3) + R_{\geq 5}$ for $p > 2$.
\end{example}

\begin{proposition}
Let $R=K[x,y,z]/(x^\delta +y^\delta +z^\delta)$ and let $\Syz(m)$
denote the sheaf of syzygies for the elements $x^3,y^3,z^3$ of total
degree $m$. Then $\Syz(5)$ is an ample sheaf for $\delta \geq 5$ and
$p \geq \delta-3$ {\rm(}or characteristic zero{\rm)}.
\end{proposition}
\proof The result follows for $\delta \geq 8$ from Corollary
\ref{amplea} and for $\delta =5$ from Example \ref{fermatfuenf}. So
suppose that $\delta =6$ or $=7$. We apply the criterion of
Hartshorne-Mumford; thus we assume that the invertible quotient
sheaf $\Syz(5) \ra \shL \ra 0$ has degree $\leq 0$. Then there
exists an invertible subsheaf of $ \Syz(5)$ of degree $\geq \delta$.
Then $\Syz(6)$ contains an invertible subsheaf $\shM$ of degree
$\geq 2 \delta $ and $\Syz(7)$ contains an invertible subsheaf
$\shM$ of degree $\geq 3 \delta $. From Riemann-Roch we see that
these subsheaves of $\Syz( \delta)$ have non-trivial global
sections.

Let $\delta =6$ and consider a global syzygy of total degree $6$. It
is given by $S_1x^3+S_2y^3+S_3z^3= \lambda F$ on the curve. The
$S_i$ do not have a common divisor, for then $\lambda = 0$ and the
syzygy would be a multiple of a syzygy of $ \Syz(k)$, $k \leq 5$,
but these are zero. This implies that the polynomials $S_1$ and
$S_2$ (say) have at most a linear form in common. Then the number of
zeros is bounded by $ 6 +4 $ or by $9$; hence it is $\leq 12 $.

Consider now $\delta =7$ and let $S_1x^3+S_2y^3+S_3z^3$ denote a
syzygy on the curve of total degree $7$; hence $\deg (S_i)=4$. The
three polynomials together have a common divisor of degree at most
$1$, since $\Gamma(Y, \Syz(5)) =0$. So two of the polynomials, say
$S_1$ and $S_2$, have a common divisor $C$ of degree $t$ at most
$2$. If $t=2$, then the number of zeros is bounded by $2 \cdot 7+
4$; if $t=1$, then it is bounded by $7 +9$; and for $t=0$ it is
bounded by $16$, so in any case $\leq 21$. \qed

\bibliographystyle{plain}

\bibliography{bibliothek}

\end{document}